\documentclass[12pt]{amsart}
\usepackage{amssymb,latexsym,amsmath,euscript}
\usepackage[T2C]{fontenc}
\usepackage{MnSymbol, mathrsfs}

\long\def\symbolfootnote[#1]#2{\begingroup%
\def\thefootnote{\fnsymbol{footnote}}\footnote[#1]{#2}\endgroup} 

\begin{document}

\title{Hypergroups and hypergroup algebras}
\author{G. L. Litvinov}
\date{}
\thanks{Journal of Soviet Mathematics vol. 38, \# 2, 1987, p. 1734-1761. \endgraf  Translated from Itogi Nauki i Tekhniki, Sovremennye Problemy Matematiki, Noveishie Dostizheniya, Vol. 26, pp. 57-106, 1985. In Russian. The translation autor is unknown.}
\maketitle

The survey contains a brief description of the ideas, constructions, results, and prospects of the theory of hypergroups and generalized translation operators.  Representations of hypergroups are considered, being treated as continuous representations of topological hyper-group algebras.

\section*{\small 1. Introduction}

{\bf 1.1.}   The important role which group-theoretic methods play in analysis and its applications, in particular in applications to theoretical physics, are well known.   Such basic mathematical concepts as translation operator, convolution, periodic function, almost periodic function, positive definite function, etc. are formulated in group-theoretic terms.   One can get far-reaching generalizations of the fundamental principles and results, connected with the concepts indicated, in the framework of the theory of hypergroups.  Essential fragments of this theory became familiar as the theory of the Delsarte -- B. M. Levitan generalized translation operators, the theory of Yu. M. Berezanskii -- S. G. Krein of hypercomplex systems with continuous basis, the theory of convolution algebras, etc.

Roughly speaking, a hypergroup is a topological space or manifold with a supplementary structure, which permits one to construct a Banach or topological algebra of the type of a group algebra -- a hypergroup algebra. Thus, in the theory of hypergroups, just as in the theory of supergroups, the object of generalization is not so much a group as a group algebra (coalgebra).

The ideas and methods of the theory of group representations carries over to the case of hypergroups, while it is convenient to treat representations of hypergroups as representations of the corresponding hypergroup algebras.   With the help of the theory of representations for hypergroups one can generalize the duality principle of L. S. Pontryagin, construct an analog of the Fourier transform, get the Plancherel theorem and the inversion formula.   It turns out that the converse result is also valid:   the existence of a transformation of the type of the Fourier transform, for which the Plancherel theorem and inversion formula are valid, is necessarily connected with the existence of a certain hypergroup.   This result explains the appearance of hypergroup structures in various problems of harmonic analysis.

For some classes of hypergroups results on a connection with infinitesimal objects in the spirit of the Lie theory are found.   Not only Lie algebras   but also algebras generated by commutation relations of a more general kind can appear as such objects.

The present survey contains a short description of the ideas, constructions, results, and prospects of the theory of hypergroups and generalized translation operators.   In composing it material of the papers [17, 57, 61, 62] has been used in part.  Separate aspects of the theory are considered in detail in the following works of monograph or survey type:   [3, 7, 37, 51, 55, 56, 66, 93, 97, 102, 103, 131, 132, 136, 160, 169, 183], in which one can find additional information and references to the literature.   An extremely large literature is devoted to hypergroups of special form -- topological groups and semigroups, and also their representations (cf. in particular, [31, 32, 40, 43, 59, 69, 72, 73, 86, 92, 156-158, 174, 177]), the systematic analysis of which would leave the framework of the present survey.   The literature cited in the present survey does not pretend to exhaustive completeness; the works included in this list contain additional bibliography.
\\

{\bf 1.2.}   The concept of hypergroup\symbolfootnote[1]{This concept was introduced independently by F. Marthy and H. S. Wall in the mid-thirties; cf. [102, 160] for the history of the question.}   arose originally as a generalization of the concept of abstract group. An abstract "algebraic" hypergroup is a set $H$ with a binary multiplication operation   $a , b \mapsto ab$   which associates with any pair of elements of $H$ a nonempty subset of $H$.   The multiplication is assumed to be associative\symbolfootnote[2]{Other definitions of abstract hypergroup occur in the literature (cf., e.g., [130, 171]); in particular, nonassociative hypergroups are considered.}  in the sense that the sets $(ab)c$ and $a(bc)$ coincide; here $(ab)c$ denotes the union of the sets $dc$ for all $d \in (ab)$, and the product $a(bc)$ is defined analogously.   A hypergroup   $H$ has an identity $e \in H$ if $a \in ea \cap ae$   for all $a \in H$.   The standard examples of hypergroups are connected with sets of cosets and conjugacy classes of elements in groups, with sets of points in certain geometries.\symbolfootnote[3]{The surveys [102, 103, 160] reflect papers devoted to the theory of abstract hypergroups.}   However, it is more convenient to start with the analysis of an example, which, at first glance, is of a different kind.

Let $G$ be a compact group,  $\widehat{G}$ be the set of all irreducible linear (finite-dimensional) representations of the group $G$, considered up to equivalence.   For any irreducible representations $\alpha$ and  $\beta$ of the group $G$, their tensor product  $\alpha \otimes \beta$ decomposes uniquely into a direct sum of primary representations
$$
\alpha \otimes \beta = \sum^n_{i=1} m_i \pi_i, \eqno(1.1)
$$
where $\pi_i \in \widehat{G}$    and $m_i$ is the multiplicity with which the irreducible representation $\pi_i$ occurs in the tensor product $\alpha \otimes \beta$.   If the product of the elements $\alpha$ and $\beta$ in $\widehat{G}$ is defined as the set $\{ \pi_1 , \pi_2 \ldots ,\pi_n \}$ of irreducible representations contained in $\alpha \otimes \beta$, then $\widehat{G}$ gets the structure of a hypergroup.

The example given was considered by Helgason in [130], which is devoted to lacunary Fourier series on compact groups.   In order to take into account the multiplicities with which the irreducible representations occur in the decomposition (1.1), Helgason defined the product $\alpha \beta$ as a finite measure on $\widehat{G}$.   The support of this measure coincides with the set $\{ \pi_1 , \pi_2 \ldots ,\pi_n \}$ of elements of $\widehat{G}$ which occur in the decomposition (1.1), and the measure of the point $\pi_i$ is the integer $m_i$.   If one identifies each element $\pi \in G$ with the unit measure $\delta \pi$, concentrated at the point $\pi$ (the delta-function), then one can consider the measure $\alpha \beta$ as the result of a convolution type operation over the measures $\delta_{\alpha}$ and $\delta_{\beta}$.   Since any measure on $\widehat{G}$ is a linear combination of delta-functions, one can extend this operation linearly to the space  $\mathscr{M} (\widehat{G})$ of all finite complex measures with finite support and even to the space   $\mathscr{M}^b (\widehat{G})$ of all bounded measures on $\widehat{G}$.   As a result, $\mathscr{M} (\widehat{G})$  and $\mathscr{M}^b (\widehat{G})$ are turned into associative algebras with identities (hypergroup algebras).
\\

{\bf 1.3.}   It is easy to modify the definition of multiplication -- convolution of measures -- in the example considered (replacing the measure $\delta_{\pi}$ by $\delta_{\pi} / \dim \pi$) so that the product of delta-functions appears as a probability measure, i.e., a positive measure with unit volume.   Constructions of this kind were studied particularly intensively after the appearance of the papers of Dunkl (112, 113], Jewett [136], Spector [169], and the beautiful survey of Ross [160].   Following these papers, in the modern literature hypergroup usually means a locally compact space $H$, such that on the set $\mathscr{M}^b (H)$  of bounded Radon measures there is given an associative bilinear operation   $\mu_1 , \mu_2 \mapsto \mu_1 \ast \mu_2$ called (generalized) convolution, where the result of the convolution of any probability measures is again a probability measure.   It is required that the convolution turn $\mathscr{M}^b (H)$ into a Banach algebra with identity, whose role is played by the delta-function $\delta_e$, concentrated at some point $e \in H$.   Moreover, it is required that the convolution be compatible with some involution in $H$ and that additional conditions of the type of continuity hold.

In particular, if $H$ is a locally compact group, then the operation in $\mathscr{M}^b (H)$ coincides with the usual convolution of measures, and $\mathscr{M}^b (H)$ coincides with the group algebra.   In this case $\delta_a \ast \delta_b = \delta_c$, where $c$ is the product of the elements $a$ and $b$ in $H$.   In general, $\delta_a \ast \delta_b$ is a probability measure which can be considered as the distribution of a random variable with values in $H$.   Hence one can say that in the hypergroup $H$ there is defined an associative product of elements, but it is defined "randomly" and its result is a "random" element in $H$.

In what follows, hypergroups in the sense of [112, 136, 169] will be called p-hypergroups, so as to distinguish them from the objects of a more general kind introduced by Delsarte.
\\

{\bf 1.4.}   Although in Ross' survey [160] it is indicated that Helgason [130] was the first analyst to use the term "hypergroup," this term was used in a broader sense in the important (and earlier) paper of Delsarte [108] in connection with various analytic applications (differential operators, mean-periodic functions, almost periodic functions).  Hypergroup in the sense of [108] is a concept which is essentially equivalent with the concept of generalized translation operators (g.t.o. for the sake of brevity) introduced by Delsarte at the end of the thirties, by the axiomatization of the properties of group translations.   Important ideas and a number of original results on the theory of generalized translation operators are due to him [105-109].   The systematic construction of the theory is given mainly in the papers of B. M. Levitan, including various questions of abstract harmonic analysis, elements of the theory of finite-dimensional linear representations, generalization of Lie theory, applications to the theory of differential operators   and questions of decomposition in eigenfunctions, generalization of the theory of almost periodic functions, etc.; cf., e.g., [47-56].

In 1950 Yu. M. Berezanskii and S. G. Krein, developing the theory of generalized translation operators (g.t.o.), introduced the concept of hypercomplex system (h.s.) with continuous basis; this basis is a locally compact space, provided with a structure which turns it into a hypergroup.   Later, harmonic analysis was constructed for h.s.   At the first stage basically commutative h.s. with compact and discrete basis were studied and the case of a general locally compact basis was to a large extent only noted; cf., e.g., [1, 2, 6-8].   The results of this stage are summarized in [7].  Recently, in [3-5, 35] a systematic construction of the theory of commutative h.s. with locally compact basis is given; [14, 15] are devoted to similar questions.   Interest in this circle of questions was reanimated under the influence of a stream of papers devoted to objects similar to  $p$-hypergroups (cf. above).   We note that in the framework of the theory of p-hypergroups, many results previously found in the theory of generalized shift operators and h.s. with continuous basis were proved again. Convolution algebras are similar to p-hypergroups and h.s. in their nature [135, 149, 150, 156, 157, 178, 20-22, 121, 122, 174].
\\

{\bf 1.5.}   In the study of hypergroups the apparatus of the theory of Banach algebras is widely applied.\symbolfootnote[1]{In the Russian literature Banach algebras are often called normed rings.}   Even in the early papers of B. M. Levitan [47-50] Banach "hypergroup" algebras are constructed; analogs of the Plancherel theorem, Bochner's theorem on the representation of positive definite functions, and the Pontryagin duality law are found for a wide class of commutative families of g.t.o.  However in certain important cases, instead of Banach algebras one should consider hypergroup algebras of a different kind.   Even for topological groups and Lie groups some class or other of topological group algebras arises naturally, depending on the version of the theory of representations and harmonic analysis being investigated, cf., e.g., [32, 59].   For hypergroups the situation is made more complex in that one is not always able to construct the Banach hypergroup algebra; a familiar example is the family of g.t.o., generated by a Sturm -- Liouville operator with rapidly increasing potential, cf. [69, Sec. 31, point 11].   Topological hypergroup algebras generated by Sturm -- Liouville operators and consisting of generalized functions (or measures) with compact support on the line were considered in [134].   A quite general theory of commutative g.t.o. in spaces of basic and generalized functions is actually constructed by L. Ehrenpreis [117].   Topological algebras of analytic functionals, connected with certain g.t.o. on the complex line are described in [41].   For g.t.o. of general type, topological hypergroup algebras and the theory of their representations are considered in [61], in which the problems of spectral analysis and synthesis are treated as problems of studying ideals of hypergroup algebras.   Hilbert hypergroup algebras are considered in [17] in connection with the generalization of the Plancherel theorem to the case of (not necessarily commutative) hypergroups.

Infinitesimal hypergroup algebras arise in generalizing Lie theory and play the same role as universal enveloping algebras of Lie algebras and associative upper-envelope Lie algebras in the sense of P. K. Rashevskii [77] for Lie groups.   B. M. Levitan proved that one can associate with each Lie algebra not only a Lie group, but also a differentiable family of g.t.o., cf. [52-56], and also [24, 25, 126]; thus enveloping algebras of Lie algebras can appear as infinitesimal objects not only for Lie groups, but also for hypergroups.  Recently, hypergroup algebras generated by non-Lie commutation relations have been investigated intensively, cf. [28-30, 36, 37, 65-67].   We note the application of the techniques of hypergroup algebras (infinitesimal and Banach), generated by commutation relations, to solving problems of mathematical physics in the framework of the operator method of V. P. Maslov [37, 65, 66].
\\

{\bf 1.6.}   A characteristic singularity of the contemporary stage of development of the theory of hypergroups is the intensive interaction of this theory with the theories of supermanifolds, supergroups, and Lie super-algebras (cf. [9, 11, 58, 64, 67, 71, 139, 162]), ringed groups (Kac algebras, cf. [38, 16, 120, 163, 164]), algebras generated by commutation relations (cf. above and [19, 33, 34, 79-81, 133, 155]), operator and topological algebras (cf. [31, 61, 62, 69, 70, 82, 120, 137, 172]), mean-periodic functions, and convolution equations (cf. [104, 116, 148, 165, 26, 27, 41, 44-46, 59, 61, 70, 75, 76, 78, 83, 87, 91, 109, 118, 119]).   The spectrum of applications of the theory of hypergroups continues to expand.   Of the relatively new domains of application of the theory we note algebraic topology (cf. [12,13, 88]) and the theory of probability, in particular, the theory of probability     measures on groups, hypergroups, and homogeneous spaces (cf. [131, 132] for a surveying account, cf. also [90, 93-95, 111, 124, 125, 139, 145, 146]).   At the same time active work continues on such traditional domains of application of the theory of generalized translation operators as differential operators and equations (cf. [7, 25, 27, 37, 511 55, 56, 65, 66, 90, 97, 98, 106, 134, 136, 145]), expansions in orthogonal systems of functions and special functions (cf. [2, 7, 51, 55, 56, 114, 131, 136, 142-144, 183]), spectral decomposition of operators (cf. [3, 5, 51, 55, 56, 135, 149, 150, 178]), duality theory (cf. [3, 7, 15, 18, 39, 50, 96, 127, 151, 152, 169]),     questions of harmonic analysis (cf. [3, 7, 14, 17, 35, 51, 55, 56, 76, 78, 87, 93, 96, 99-101, 107, 136, 140-144, 167-170, 175, 177, 179-183]), etc.
\\

{\bf 1.7.}   The author considers it his pleasant duty to thank V. M. Bukhshtaber, D. I. Gurevich, B. M. Levitan, L. I. Vainerman, A. M. Vershik, and D. P. Zhelobenko  for valuable suggestions, assistance, and guidance.

\section*{\small 2.     Hypergroups   and   Generalized   Translation  Operators}

\underline{{\bf 2.1.}   Basic Concepts.}   Let $H$ be an arbitrary set, $\Phi$ be some linear space of complex-valued functions, defined on $H$.   Let us assume that to each element $x \in H$ there is assigned a linear operator $R^x$ in $\Phi$, where for any fixed $y \in H$ the function $\psi (x) =R^x \varphi (y)$ is contained in $\Phi$ for all $\varphi \in \Phi$.   We denote the linear operator $\varphi (x) \to \psi (x) = R^x \varphi (y)$ in $\Phi$ by $L^x$.

The set $H$ is called a {\it hypergroup}, if the following conditions (axioms) hold:

1)   Associativity Axiom. For any elements $x, y \in H$ one has
$$
R^x L^y = L^y R^x .
$$ 
2)   There exists in $H$ an element $e$, called neutral (or the identity), such that $R^e = I$, where $I$ is the identity operator in $\Phi$.

In this case one says that the operators $R^x$ form a family of generalized translation operators (g.t.o.); thus the concepts of g.t.o. and hypergroup are considered as equivalent.   The operators $R^x$ are often called (generalized) right translation operators, in which case the $L^y$ are called left translation operators.

G.t.o. obviously arise in any linear subspace of functions on an arbitrary group or semigroup with identity which is invariant with respect to translations.   We set $R^x \colon \varphi (h) \mapsto \varphi (hx)$ where $hx$ is the product of the elements $h$ and $x$ in the semigroup.   It is clear that $L^y \varphi (h) = \varphi (yh)$ and the axiom of associativity reduces to the associativity of multiplication in the semigroup, and the neutral element is the identity in the semigroup.  In the present case the operators $L^y$ also form a family of g.t.o.

However in general the $L^y$ do not form g.t.o., since the operator $L^e$ is not necessarily the identity.   It is easy to verify that this operator is a projector; its image $\widetilde{\Phi}$ is called the {\it basic subspace} of $\Phi$.  A basic subspace is invariant with respect to left and right translations, and in $\widetilde{\Phi}$ the operators $L^y$ form a family of g.t.o., so that the symmetry between left and right translations is reestablished.   Often the second axiom is strengthened by requiring that $L^e = I$, i.e., $\Phi = \widetilde{\Phi}$, or weakened, for example, by requiring the existence of only an approximate identity (in some sense or other).

Conditions 1) and 2) are the most general axioms for a hypergroup.   By adding additional conditions, one can single out narrower classes of hypergroups.   For example, if
$$
R^x R^y = R^y R^x \eqno (2.1)
$$
for all $x, y \in H$, then the hypergroup $H$ is said to be {\it commutative} (respectively, the g.t.o. $R^x$ are called {\it commutative}).   In this case $R^x \varphi = L^x \varphi$ for any function $\varphi$ from the basic subspace $\widetilde{\Phi}$.

If $H$ is a locally compact space\symbolfootnote[1]{In what follows all locally compact spaces and manifolds will be assumed (without saying anything special) Hausdorff and countable at infinity, i.e., admitting a representation as a union of a no more than countable family of compact subsets.} with a positive measure $m$, then it is usually required that the operators $L^x$ and $R^x$ act compatibly on the space $C(H)$ of continuous functions on $H$ and on the spaces $L^p (H,m)$ for $p \geq 1$, while on $R^x$ and $L^x$ one imposes additional conditions of the type of continuity; if $H$ is a smooth manifold, then one imposes conditions of the type of differentiability, etc.   Cf. [3, 7, 14, 17, 18, 47, 51, 55, 56, 61, 65, 105, 108, 112, 113, 136, 160, 169, 156, 157] for various versions of the axiomatics of hypergroups and g.t.o.
\\

\underline{{\bf 2.2.}   Example.  Hypercomplex Systems (h.s.).}   Let $\Phi$ be a finite h.s., i.e., a finite-dimensional associative algebra with fixed basis $H = \{ h_1 ,\\ \ldots , h_n \}$.   We identify $\Phi$ with the space of functions on the finite set $H$, by assigning to the function $\varphi$ the element with coordinates $\varphi (h_1) , \ldots ,\\ \varphi(h_n)$, i.e., the element $\sum^n_{i=1} \varphi (h_i) h_i \in \Phi$. We set $R^x \varphi (h) = \sum^n_{i=1} \varphi (h_i) h_i \ast x$ where $h_i \ast x$ is the product of the elements $h$ and $x$ in the algebra $\Phi$.   It is easy to verify that the operators $R^x$ form a family of g.t.o. (so that $H$ is a hypergroup) if and only if one of the elements of the basis of $H$ is a right identity in the algebra $\Phi$ (a two-sided identity if one requires that the basic subspace $\widetilde{\Phi}$ coincide with $\Phi$).   In the way indicated one establishes a correspondence between finite hypergroups and finite h.s.

Thus, the concept of hypergroup (as well as that of g.t.o.) can be considered as a far-reaching generalization of the classical concept of hypercomplex system.   Many examples of infinite hypergroups which are naturally treated as h.s. with countable or continuous basis are considered in [1-8, 51, 55, 56]  and in the papers devoted to p-hypergroups (cf. the Introduction).   Some of these examples are considered below.
\\

\underline{{\bf 2.3.} Delsarte's Example (cf. [108]).}   Let $G$ be a topological group, $K$ be a compact group of automorphisms of the group $G, dk$ be an invariant measure on $K$, and $\int dk = 1$.   In the space $\Phi = C(G)$ of continuous functions on $G$, we define a g.t.o. $R^x$ with the help of the equation:
$$
R^x \varphi (g) = \int_K \varphi (g \cdot k (x)) dk, \eqno(2.2)
$$
where $\varphi \in \Phi, k(x)$ is the image of the element $x \in G$ under the automorphism $k \in K, g \cdot k(x)$ is the product of the elements $g$ and $k(x)$ in $G$. It is easy to verify that both hypergroup axioms hold, where the neutral element is the identity of the group $G$.   The hypergroup constructed is commutative, if the group $G$ is commutative.

In the present case the basic subspace $\widetilde{\Phi}$ consists of all functions which are constant on orbits with respect to the action of the group $K$, and the operators $R^x$ and $R^y$ coincide in $\widetilde{\Phi}$, if $x$ and $y$ lie on the same orbit. Hence the space of orbits $H$ can also carry the structure of a hypergroup, identifying $\widetilde{\Phi}$  with the space $C(H)$ of continuous functions on $H$ (with respect to the quotient toplogy) and setting $R^h = R^x$, where $x$ is an arbitrary element of the orbit $h$.

The process described, of constructing a hypergroup consisting of orbits, is called reduction.  In general, reduction consists of replacing $\Phi$ by $\widetilde{\Phi}$ and identifying elements of the original hypergroup if the values of any function on $\widetilde{\Phi}$ on these elements coincide.   For the reduced hypergroup one has a strengthened form of the axiom of existence of an identity: $R^e = L^e = I$.   On the other hand, as a result of reduction, a smooth manifold can be turned, for example, into a manifold with boundary (cf. below point 2.4).

A generalization of Delsarte's construction, which lets one get a large class of hypergroups, is given in the important paper [28].
\\

\underline{{\bf 2.4.}   Example.}   We consider a special case of the construction described in point 2.3, when $G$ coincides with the real axis ${\bf R}$, and the group of automorphisms consists of two elements (reflection with respect to zero and the identity map).   In this case (2.2) takes the form
$$
R^x \varphi (t) = \frac{1}{2} [\varphi (t + x)+ \varphi (t-x)] . \eqno(2.3)
$$
$x, t \in {\bf R}$.   The basic subspace consists of all even functions.   The operator of left generalized translation has the form
$$
L^y \varphi (t) = \frac{1}{2} [ \varphi (t+y)+ \varphi (-t+y)] .
$$
It is clear that the operator $L^y$ for any value of the parameter $y$ carries an arbitrary function into an even one. It is easy to verify that the given hypergroup is commutative and on the basic subspace the right generalized translations coincide with the left ones.   Reduction also provides the orbit space, which can be identified with the half-line $0 \leq t < \infty$.
\\

\underline{{\bf 2.5.}   Example.}   Another interesting special case of Delsarte's example we have for $G = K$ and $k(x) = k^{-1} xk$; here the basic subspace consists of all central functions on $G$, i.e., functions which are constant on conjugacy classes of elements.   The reduced hypergroup $H$ is commutative and is made up of conjugacy classes of elements of the compact group $G$.   This hypergroup is closely connected with representations of the group $G$ and is dual with respect to the commutative hypergroup $\widehat{G}$, considered above in point 1.2.   The hypergroup algebra for $H$ coincides with the center of the group algebra $G$.   From this point of view the hypergroup of conjugacy classes of elements of a compact group was studied in [1].   Cf. [160, 161, 129] for various generalizations.   For example, if a compact group of automorphisms $K$ of an arbitrary locally compact group $G$ includes all inner  automorphisms, then the orbit space in $G$ with respect to $K$ is provided with a commutative hypergroup structure.
\\

\underline{{\bf 2.6.}   Double Cosets with respect to a Compact Subgroup.}   Let $G$ be a locally compact group with a compact subgroup $K$, and $H = K\setminus G/K$   be the space of double cosets with respect to the subgroup $K$ (such a coset contains along with the element $g \in G$ all elements of the form $k_1 gk_2$, where $k_1 , k_2 \in K$).  If $K$ is a normal subgroup of $G$, then $H$ coincides with the quotient group $G/K$.   Let $\Phi$ be the space consisting of all those continuous functions $\varphi (g)$ on $G$, such that $\varphi (k_1 g k_2) = \varphi (g)$ for any elements $k_1 , k_2 \in K, g \in G$.   The family $R^x$ of right generalized translations in $\Phi$ is defined by
$$
R^x \varphi (g) = \int_K \varphi (gkx) dk . \eqno(2.4)
$$
The space $\Phi$ can be identified with the space $C(H)$ of continuous functions on $H$, and arguing just as in Delsarte's example (cf. above point 2.3), one can provide $H$ with a hypergroup structure.   We note that (2.4) defines a g.t.o. in the space of continuous functions on $G$, invariant with respect to the action of the subgroup $K$ on the right, also, i.e., on the space of continuous functions on the homogeneous space $G/K$ of left cosets.   The hypergroup $H$ is constructed from the hypergroup $G/K$ by reduction, while elements in $H$ correspond to orbits in $G/K$ with respect to the action of the stationary subgroup $K$.

For example, if $G$ is the group of Euclidean motions of the plane, and $K$ is the subgroup of rotations about a fixed point, then $G/K$ coincides with the plane, the orbits in $G/K$ with respect to the action of $K$ form a family of concentric circles and an element of the hypergroup $H = K\smallsetminus G/K$ is determined by the radius of the corresponding circle.   Thus, in the present case H can be identified with the half-line $[0, \infty )$.   This hypergroup is closely connected with the Helmholtz equation on the plane and the Bessel equation (cf., e.g., [136]).  Starting from different considerations, this same example was actually already considered by Delsarte [105].

Hypergroups of this type were studied in the frameworks of the theory of hypercomplex systems with continuous basis and the theory of p-hypergroups (cf., e.g., [7, 131, 136, 169, 146]).
\\

\underline{{\bf 2.7.}  Adjoints of a Generalized Translation Operator (g.t.o.) and}\\ \underline{Hypergroups with Involution.}   Let the hypergroup $H$ be a locally compact space, on which there is given a positive Borel measure $m$.   Let us assume that the action of the g.t.o. $R^x$ and $L^x$ is defined not only on the original space $\Phi$, but also on the Hilbert space $L^2 (H, m)$.  We denote by $\tilde{R^x}$ the adjoint of the operator $R^x$, defined by
$$
\int_H R^x \varphi (y) \cdot \bar{\psi} (y) dm(y) = \int_H \varphi (y) \cdot \overline{\widetilde{R^x} \psi (y)} dx (y) , \eqno(2.5)
$$
where $\varphi, \psi \in L^2 (H, m)$ and $\psi \mapsto \bar{\psi}$  is complex conjugation.   The operators $\widetilde{R^x}$ form a family of adjoint (right) g.t.o.; one defines adjoint operators of left translation $\widetilde{L^x}$ analogously.   If $G$ is a locally compact group with right invariant measure $m$, then the adjoint of the right translation $\varphi (g) \mapsto \varphi (gx)$ is the operator $\varphi (g) \mapsto \varphi (g x^{-1})$; if the measure m is left-invariant, then the operator adjoint to left translation by $x$ has the form $\varphi (g) \mapsto \varphi (x^{-1} g)$.

B. M. Levitan formulated a condition on the measure $m$ of the type of invariance of this measure with respect to generalized translations in terms of adjoints of g.t.o. (cf. [47, 51, 55, 56]).   In particular, for commutative g.t.o. the corresponding condition has the form
$$
\widetilde{R^x} R^y = R^y \widetilde{R^x} , \eqno(2.6)
$$
i.e., the g.t.o. and the adjoint g.t.o commute; if (2.6) holds, then the hypergroup is called normal.

Let us assume that there exists an involutive homeomorphism $x \mapsto x^*$ mapping the hypergroup $H$ to itself, while in $\Phi$ there is induced an involution $\varphi (x) \mapsto \tilde{\varphi} (x) = \varphi (x^*)$. We call $H$ a hypergroup with involution, if for all $x, y \in H, \varphi \in \Phi$ one has
$$
(R^{x^*} \varphi ) (y^*) = \overline{R^y (\Tilde{\Tilde \varphi} (x))} . \eqno(2.7)
$$
If the g.t.o. $R^x$ are real, i.e., commute with complex conjugation, then (2.7) assumes the form
$$
(R^{x^*} \varphi ) (y^*) = R^y (\tilde \varphi (x)) . \eqno(2.8)
$$
(2.7) and (2.8) mean that an involution is also induced on the hypergroup algebra (cf. below, point 3.6). An involution on a hypergroup is a substitute for the map $x \mapsto x^{-1}$ in a group. For a hypergroup with involution the space $\Phi$ coincides with its basic subspace $\widetilde{\Phi}$.

We call a measure $m$ on $H$ right-invariant, if the adjoint right translation operator $\widetilde{R^x}$ coincides with the right translation $R^{x*}$ , defined by the element $x^* \in H$; similarly a measure is left-invariant, if $\widetilde{L^x} = L^{x*}$.   Under additional conditions, which reduce to the fact that the function $\varepsilon (x) \equiv 1$ (as well as other constants) is unchanged under the action of the g.t.o., a right-invariant (left-invariant) measure is unchanged under the action of right (left) translations.   In fact, formally substituting the function $\varepsilon$ into (2.5) in place of $\psi$ and considering that $\widetilde{R^{x*}} \varepsilon = R^{x*} \varepsilon \equiv 1$ we get $\int R^x \varphi (y) dm (y) = \int \varphi (y) dm (y)$ which is what was required.  With the help of the involution a left-invariant measure can be transformed into a right-invariant one and conversely.

The reduced hypergroups described above in points 2.3 -- 2.6 are hypergroups with involutions.   If $g$ is an element of the orbit (coset) $x$, then $x^*$ coincides with the orbit (coset) of the element $g^{-1}$.   The g.t.o. arising here are real.   The hypergroups indicated have invariant measures, induced by the invariant measures on the corresponding groups.

A hypergroup with involution is called {\it Hermitian}, if the involution reduces to the identity map. For example, the reduced hypergroup described in point 2.4 is Hermitian. A Hermitian hypergroup is automatically commutative, if it generates real g.t.o.

\section*{\small 3. Hypergroup   Algebras}

{\bf 3.1.}   If V is a locally convex space\symbolfootnote[2]{To be specific, we shall consider all linear spaces and algebras to be over the field of complex numbers. Nevertheless many of the results described below are also valid for spaces and algebras over the field of real numbers.},  and $V'$ is the dual space (consisting of all continuous linear functionals on $V$), then the value of the functional $v' \in V'$ on the element $v \in R$ will be denoted by $\langle v', v \rangle$.   We recall that the weak topology $\sigma (V', V)$ on $V'$ is defined by the seminorms $\upsilon ' \mapsto \vert \langle \upsilon ', \upsilon \rangle \vert $ , where $v$ runs through $V$; analogously, the seminorms  $\upsilon ' \mapsto \vert \langle \upsilon ', \upsilon \rangle \vert $, where $v' \in V'$, define the weak topology on $V$.   The strongest locally convex topology on $V$, admitting $V'$ as the space of all continuous linear functionals, is called the {\it Mackey topology} and is denoted by $\tau (V, V')$; one defines the {\it Mackey topology} $\tau (V', V)$ on $V'$ analogously. The space $V$ is called a {\it Mackey space} if its topology is the Mackey topology $\tau (V, V')$. All the barreled spaces (in particular, all the Banach spaces and all the complete metrizable spaces) are Mackey spaces.  In what follows, if nothing is said to the contrary, then dual spaces are provided with the Mackey topologies.   If $V$ is a Mackey space, then passage to the Mackey-dual space is reflexive in the sense that $(V')' = V$.   We denote by $S(V)$ the algebra of all weakly continuous linear operators on $V$, provided with the weak operator topology; this topology is defined by the seminorms $A \mapsto \vert \langle \upsilon ', A \upsilon \rangle \vert$, where $v$ and $v'$ run through $V$ and $V'$ respectively.   If $V$ is a Mackey space, then any operator from $S(V)$ is continuous.   Elements of $S(V)$ are denoted below either by capital Latin or by lower-case Greek letters.

Let $\Phi$ be a locally convex space, consisting of functions $\varphi (x)$, defined on some set $H$.   The value $\langle f, \varphi \rangle$ of the functional $f \in \Phi '$ on the function $\varphi \in \Phi$ will sometimes be denoted by $\int \varphi (x) d f (x)$.  If the operator-valued function $x \mapsto A^x$ maps $H$ into $S(V)$ and the function $\varphi (x) = \langle v', A^x v \rangle$ is contained in $\Phi$ for all $v \in V,  v' \in V'$, then we denote by $A (f) = \int A^x df (x)$   the linear operator on $V$, defined by
$$
\langle \upsilon ', A(f) \upsilon \rangle = \int \langle \upsilon ', A^x \upsilon \rangle df(x) , \eqno(3.1)
$$
where $v$ and $v'$ run respectively through $V$ and $V'$ (the existence of the operator $A(f)$ in each concrete case requires a proof).   The notation $\int \pi_x df(x)$ has analogous meaning, where $x \mapsto \pi_x$ is an operator-valued function.

The information on locally convex spaces needed for what follows can be found, for example, in [89].   The proofs of the basic assertions on hypergroup algebras formulated in points 3.2 -- 3.6 are easily extracted from [61, 62]; cf. also [156, 157].
\\

\underline{{\bf 3.2.}   Measure Algebras on Locally Compact Hypergroups.}   We consider g.t.o. acting on the space $C(H)$ of continuous functions on the locally compact space $H$, where $C(H)$ is provided with the topology of uniform convergence on compacta\symbolfootnote[3]{Since $H$ is assumed to be countable at infinity (cf. the footnote on p. 1737), $C(H)$ is a complete metrizable space and consequently a Mackey space.}.  Let us assume that the right translation and left translation operators $R^x$ and $L^x$ are continuous and depend weakly continuously on the parameter $x$; in this case we shall say that the given family of g.t.o. is {\it continuous} and generates in $H$ a {\it locally compact hypergroup structure}.   For this it is sufficient that the correspondence $\varphi (x) \mapsto \psi (x,y) = R^x \varphi (y)$ be a continuous map $C(H) \to C(H \times H)$.

Let $\mathscr{M} (H)$ be the dual space to $C(H)$, consisting of complex Radon measures with compact support. For any element $f \in \mathscr{M} (H)$, with the help of equations of the type of (3.1), one defines continuous linear operators on $C(H)$:
$$
R(f) = \int R^x df(x), L(f) = \int L^x df(x) . \eqno(3.2)
$$
By the generalized convolution (or simply convolution) of the elements $f, g \in \mathscr{M} (H)$ is meant the functional $f \ast g$ defined by
$$
\langle f \ast g, \varphi \rangle = \langle f, R(g) \varphi \rangle . \eqno(3.3)
$$
Since the operator $R(g)$ is continuous, $f \ast g \in \mathscr{M} (H)$. One can verify that one also has
$$
\langle f \ast g, \varphi \rangle = \langle g, L(f) \varphi \rangle . \eqno(3.4)
$$
(3.3) and (3.4) show that the operator $R(g)$ is adjoint\symbolfootnote[2]{We recall that the adjoint operator $A' \colon V' \to V'$ of the (at least weakly) continuous linear operator $A$ on the locally convex space $V$ is defined by the equation $\langle A'v', v \rangle = \langle v', Av \rangle$ for all $v \in V, v' \in V'$.   The adjoint operator $A'$ depends linearly on $A$.   Above, in point 2.7, we considered the transition $A \to \widetilde{A}$ to the adjoint operator on a Hilbert space.   One should keep in mind that this operation is antilinear, since the Hilbert space is identified with its dual by means of an antilinear correspondence.} to the operator $f \mapsto f \ast g$ on   $\mathscr{M} (H)$ and the operator $L(f)$ is adjoint to the operator $g \mapsto f \ast g$. It follows from this and from the continuity of the operators $L(f), R(g)$, that convolution is separately continuous in the weak topology $\sigma (\mathscr{M} (H), C(H))$ and in the Mackey topology $\tau (\mathscr{M} (H), C(H))$, i.e., gives a separately continuous bilinear map $\mathscr{M} (H) \times \mathscr{M} (H) \to \mathscr{M} (H)$.   It is easy to prove that the operators $L(f)$ and $R(g)$ commute for all $f, g \in \mathscr{M} (H)$   from which it follows that convolution is an associative operation.   Thus, $\mathscr{M} (H)$ is a topological algebra (or locally convex algebra), i.e., an associative algebra, provided with a locally convex topology such that the product of elements is separately continuous.

We denote by $\delta_x$ the Dirac measure, concentrated at the point $x \in H$ ("$\delta$ -function"), i.e., the functional from $\mathscr{M} (H)$ such that $\langle \delta_x , \varphi \rangle = \varphi (x)$ for all $\varphi \in C(H)$.    It is well known that the element  $f \in \mathscr{M} (H)$ coincides with some Dirac measure $\delta_x$,if and only if it is a multiplicative functional on $C(H)$, i.e., when $\langle f, \varphi  \psi \rangle = \langle f, \varphi \rangle \langle f. \psi \rangle$ where $\varphi  \psi$ is the usual product of functions $\varphi , \psi \in C(H)$.

It is clear from (3.2) that $R(\delta_x ) = R^x$   and $L( \delta_x ) = L^x$.    If $e$ is the neutral element of $H$, then $\langle f \ast \delta_e , \varphi \rangle = \langle f, R^e \varphi \rangle = \langle f, \varphi \rangle$ for all $f \in \mathscr{M} (H) , \varphi \in C(H)$, i.e., $f \ast \delta_e = f$  and $\delta_e$ is a right identity in $\mathscr{M} (H)$.

Thus, the following result is valid.
\\

\underline{THEOREM 1.} For any locally compact hypergroup $H$ the operation of generalized convolution turns the space $\mathscr{M} (H)$ provided with the Mackey topology (or the weak topology) into a topological algebra with right identity, where this identity is a multiplicative functional on $C(H)$.
\\

The converse theorem is also valid.
\\

\underline{THEOREM 2.}   Let $H$ be a locally compact space and suppose given in $\mathscr{M} (H)$ the structure of a topological algebra with right identity, while this identity is a multiplicative functional on $C(H)$, and the topology in $\mathscr{M} (H)$ is compatible with the duality between $\mathscr{M} (H)$ and $C(H)$ (for example, coincides with the Mackey or with the weak topology).   Then there exists a unique locally compact hypergroup structure on $H$ such that the multiplication in $\mathscr{M} (H)$ coincides with the corresponding convolution.
\\

If there is given in $\mathscr{M} (H)$ the structure of a topological algebra in accord with the hypothesis of Theorem 2, then one can define the operator of right translation $R^x$ in $G(H)$ as the adjoint of the operator of multiplication by $\delta_x$ on the right in $\mathscr{M} (H)$.   Respectively, $L^x$ is adjoint to the operator of multiplication by $\delta_x$ on the left.   The adjoint to the projector $L^e$ to the basic subspace $\widetilde{C} (H)$ in $C(H)$ (cf. above, point 2.1)   is the projector $\tau \mapsto \delta_e \ast f$ in $\mathscr{M} (H)$. The image of this operator will be denoted by $\widetilde{\mathscr{M}} (H)$.   It is easy to prove
\\

\underline{Proposition 1.} $\widetilde{\mathscr{M}} (H)$ is a closed subalgebra in  $\mathscr{M} (H)$ while $\delta_e$ is a two-sided identity in $\widetilde{\mathscr{M}} (H)$. One can identify the subalgebra $\widetilde{\mathscr{M}} (H)$ with the quotient-algebra (provided with the quotient-topology) of the algebra $\mathscr{M} (H)$ by the ideal orthogonal to the basic subspace $\widetilde{C} (H)$.

Thus, $\widetilde{\mathscr{M}} (H)$ coincides with the dual space to $\widetilde{C}(H)$.   The algebra $\widetilde{\mathscr{M}} (H)$ is called the basic subalgebra of $\mathscr{M} (H)$.    This algebra is the hypergroup algebra of the reduced hypergroup (cf. above, point 2.3).  We note that the hypergroup algebra $\mathscr{M} (H)$ is noncommutative if $\widetilde{C} (H)$ does not coincide with $C(H)$, even when the hypergroup $H$ is commutative.   On the other hand, in this case the basic subalgebra $\widetilde{\mathscr{M}} (H)$ is commutative.

The constructions and results on hypergroup algebras described can be considered as the formalization of the heuristic considerations expressed by B.M. Levitan (cf., e.g., [54, Sec. 2, point 1]).   These constructions and results can also be carried over to the case when instead of $C(H)$ one considers the Banach space $C_0 (H)$, consisting of all continuous functions on $H$, which vanish at infinity.  In other words, $\varphi (x) \in C_0 (H)$ if for any $\varepsilon > 0$ one can find a compact set $K \subset H$ such that $\vert \varphi (x) \vert < \varepsilon$  for $x \notin K$ and the norm is determined by the equation $\| \varphi \| = \max_{x \in H} | \varphi (x)|$.     The dual space $\mathscr{M}^{b} (H)$ of $C_0 (H)$ consists of all bounded complex Radon measures.
If in $C_0 (H)$ there are defined continuous left and right generalized translation operators $L^x$ and $R^x$, where the functions $x \mapsto \langle f, R^x \varphi \rangle$ and $x \mapsto \langle f, L^x \varphi \rangle$ are contained in $C_0 (H)$ for all $F \in \mathscr{M}^{b} (H), \varphi \in C_0 (H)$, then $\mathscr{M}^{b} (H)$ can be provided with the structure of a topological hypergroup algebra and even a Banach algebra (with respect to to the norm of the Banach space dual of $G_0 (H)$).

Structures of the type indicated on $\mathscr{M}^{b} (H)$ were apparently first considered for a sufficiently general case in [156, 157] in connection with the theories of semigroups and semigroup algebras.   The algebra $\mathscr{M}^{b} (H)$ is a basic object in the theory of $p$-hypergroups (cf., e.g., [112, 113, 136, 169, 160, 93, 131, 132]).
\\

\underline{{\bf 3.3.}   Examples.}   In Delsarte's example, described in point 2.3, the hypergroup algebra $\mathscr{M} (H)$ coincides with the subalgebra of the group algebra $\mathscr{M} (G)$ consisting of measures, invariant with respect to the action of the compact group $K$.   For the special case of Delsarte's example described in point 2.4, $\widetilde{\mathscr{M}} (R)$  consists of all even measures, while for these measures the generalized convolution coincides with the ordinary one. For the hypergroup of classes of conjugate elements, described in point 2.5, the hypergroup algebra consists of central measures on $G$.   Finally, if $H$ is the hypergroup of double cosets, described in point 2.6, then $\mathscr{M} (H)$ is the subalgebra of the group algebra  $\mathscr{M} (G)$ consisting of measures, invariant with respect to left and right translations by elements of the compact subgroup $K$.   The hypergroup algebras of bounded measures admit an analogous description in all the examples considered.
\\

\underline{{\bf 3.4.}   Differentiable and Analytic Hypergroups.}   Let $H$ be an infinitely differentiable manifold, $C^{\infty} (H)$ be the space of all infinitely differentiable functions on $H$, provided with the topology of uniform convergence along with derivatives on compact subsets of $H$.   The dual space $\mathscr{D} (H)$ of $C^{\infty} (H)$ consists of all generalized functions with compact support on $H$.   The space $C^{\infty} (H)$ is reflexive (and is a Mackey space); hence the Mackey topology on $\mathscr{D} (H)$ coincides with the strong topology of uniform convergence on bounded sets from $C^{\infty} (H)$. A family of g.t.o. $R^x$ on $C^{\infty} (H)$ is called infinitely differentiable (respectively, $H$ is called an infinitely differentiable hypergroup), if the correspondence $\varphi (x) \mapsto \psi (x,y) = R^y \varphi (x)$ is a continuous map  $C^{\infty} (H) \to C^{\infty} (H \times H)$. It follows from this that the operators $R^x$ and $L^x$ are continuous, and the functions $x \mapsto \langle f, R^x \varphi \rangle$    and $x \mapsto \langle f, L^x \varphi \rangle$ are infinitely differentiable for any $ \varphi \in C^{\infty} (H), f \in \mathscr{D} (H)$.   Arguing as in point 2.2, one can define a generalized convolution in $\mathscr{D} (H)$.    Strengthened analogs of Theorems 1 and 2 are valid here.
\\

\underline{THEOREM 1'.}   The generalized convolution turns  $\mathscr{D} (H)$ into an associative topological algebra with continuous (and not only separately continuous) multiplication and with a right identity which is a multiplicative functional on $C^{\infty} (H)$.
\\

\underline{THEOREM 2'.}   Suppose given in $\mathscr{D} (H)$ an associative algebra structure and a multiplicative functional which is a right identity in $\mathscr{D} (H)$.    If the multiplication in $\mathscr{D} (H)$ is separately continuous in any locally convex topology, which is compatible with the duality between $\mathscr{D} (H)$ and $C^{\infty} (H)$, then this multiplication is also continuous in the Mackey topology   $\tau (\mathscr{D} (H), C^{\infty} (H))$, which coincides with the strong topology of the dual space to $C^{\infty} (H)$.   In this case there exists a unique structure of infinitely differentiable hypergroup on $H$, such that the multiplication in $\mathscr{D} (H)$ coincides with the corresponding convolution.
\\

The analog of Proposition 1 is also valid, so that one can define the basic subalgebra of  $\mathscr{D} (H)$ and it is natural to consider this subalgebra as the hypergroup algebra of the reduced hypergroup.

Let us now assume that $H$ is a complex analytic manifold.   We denote by $\mathscr{H} (H)$ the space of all holomorphic functions on $H$ with the topology of uniform convergence on compact subsets, and by $\mathscr{A} (H)$ the dual space, consisting of analytic functionals.   In $H$ there is given a structure of {\it holomorphic (or complex-analytic) hypergroup} if on $\mathscr{H} (H)$ there act g.t.o. $R^x$ such that the map $\varphi (x) \mapsto \psi (x,y) = R^y \varphi (x)$ is a continuous map $\mathscr{H} (H) \mapsto \mathscr{H} (H \times H)$.   In this case $\mathscr{A} (H)$  is turned into a hypergroup algebra with respect to the generalized convolution and the analogs of Theorems 1 and 2 and Proposition 1 are valid.

Finally, if the g.t.o. act on the space of real analytic functions on a real analytic manifold, then (if the standard conditions hold) one can endow the space of hyperfunctions with compact support with the structure of a topological hypergroup algebra.
\\

\underline{{\bf 3.5.}  Dual Topological Algebras and Topological Hopf Algebras.}   The results formulated above admit further formalization which clarifies their nature.   Let $\Phi$ be a topological algebra (associative and with separately continuous map) and let the dual space $\mathscr{F}$ of $\Phi$  also be endowed with the structure of a topological algebra, where $\Phi$ and $\mathscr{F}$ are given the Mackey topologies $\tau (\Phi , \mathscr{F} )$ and $\tau (\mathscr{F} , \Phi )$.    We call the pair $\mathfrak{A} = (\Phi , \mathscr{F} )$ a {\it double topological algebra} (for short d.t.a.).   Wanting to indicate that the algebra $\Phi$ has a certain property (for example, that it is commutative or has an identity), we shall say that the d.t.a. $\mathfrak{A}$   has this property; now if we want to indicate that the algebra $\mathscr{F}$ has some property, then we shall use the prefix "co."   For example, if $H$ is a locally compact hypergroup, then the pair $\mathfrak{A} = (C(H), \mathscr{M} (H))$ where the multiplication in $C(H)$ coincides with the multiplication of functions, and in $\mathscr{M} (H)$ with the generalized convolution, forms a commutative d.t.a. with identity and right coidentity.   Analogously, one defines d.t.a. connected with infinitely differentiable and holomorphic hypergroups.   In the theory of hypergroups, d.t.a. play the same kind of role as double Hilbert algebras do in the theory of ringed groups (Kac algebras), cf. [38, 16, 120, 163, 164].

If the algebra $\Phi$ is commutative, then we denote by $H = H(\Phi )$ the subset of   $\mathscr{F}$ consisting of all multiplicative functionals on $\Phi$ and endowed with the induced topology.   The Gelfand transformation $\varphi \mapsto \varphi (x) = \langle x, \varphi \rangle$ maps $\Phi$ into the space $C(H)$ of continuous functions on $H$.   We call the algebra $\Phi$ semisimple, if the Gelfand transformation has no kernel\symbolfootnote[1]{Usually semisimple algebras are defined differently, but for commutative multiplicative algebras which are convex in the sense of [153] (in particular, for Banach algebras), the definition given is equivalent with the ordinary one.}; in this case we identify $\Phi$ with a linear subspace of $C(H)$.  If $x \in H$, then by $R^x$ (respectively, $L^x$) we denote the operator on $\Phi$, adjoint to multiplication by $x$ on the right (left) in $\mathscr{F}$.
\\

\underline{Remark (2011).} The concept of d.t.a. is closely related with the modern concept of quantum hypergroup (and quantum group).
\\

\underline{Proposition 2.} If the d.t.a. $\mathfrak{A}$  is commutative, semisimple, and has a right coidentity, which is a multiplicative functional, then the $R^x$ are g.t.o. and multiplication in $\mathscr{F}$ coincides with generalized convolution; here the $L^x$ are left translation operators.   Thus $H$ is endowed with the structure of a hypergroup, and $\mathscr{F}$ is a hypergroup algebra.
\\

\underline{Proposition 3.}   If the condition of Proposition 2 holds, then all the operators $R^x$ are endomorphisms of the algebra $\Phi$ if and only if the multiplication in  $\mathscr{F}$   induces a semigroup structure on $H$; here $R^y \varphi (x) = \varphi (xy)$. The semigroup indicated is a group, if $\mathfrak{A}$ has a two-sided coidentity $e \in H$ and there exists an involutive automorphism $S \colon \Phi \to \Phi$, such that $\langle x, S R^x \varphi \rangle = \langle e, \varphi \rangle$; in this case $(S \varphi ) (x) = \varphi (x^{-1})$.   The group $H$ is topological (i.e., the multiplication in $H$ is not only separately continuous, but is also continuous), if the space $H$ is locally compact.
\\

We call the d.t.a. $\mathfrak{A}' = (\mathscr{F} , \Phi )$   dual to the d.t.a. $\mathfrak{A} = (\Phi , \mathscr{F} )$.    It is clear that  $(\mathfrak{A}' )' = \mathfrak{A}$.    If the d.t.a. $\mathfrak{A}$  is commutative and cocommutative, semisimple and cosemisimple, has an identity and a coidentity, which are multiplicative functionals\symbolfootnote[2]{We recall that in point 2.7 the condition that the function $\varepsilon (x) \equiv l \in \Phi$ is invariant with respect to translations was considered.   This condition means that the identity in $\Phi$ is a multiplicative functional on $\mathscr{F}$.}, then corresponding to Proposition 2 the d.t.a. $\mathfrak{A}$  generates a hypergroup $H$, and the d.t.a. $\mathfrak{A}'$  generates a hypergroup $H'$, which is naturally considered as dual to $H$.   The hypergroups $H$ and $H'$ are commutative and $(H')' = H$.   The present version of the duality theory is not a direct generalization of Pontryagin duality\symbolfootnote[3]{Questions connected with duality are touched on only episodically in the present paper.   A survey of L. I. Vainerman, "Duality for algebras with involution and generalized translation operators," will be devoted to this subject in Mathematical Analysis, Vol. 24, Itogi Nauki i Tekhniki VINITI AN SSSR, Moscow (1986).}.   For example, if  $\mathfrak{A} = (C^{\infty} ({\bf R}), \mathscr{D} ({\bf R}))$   where $\mathscr{D} ({\bf R})$  is considered as the group algebra of the group $R$, then the differentiable hypergroup $H$ coincides with this group, and the dual hypergroup $H'$ coincides with the group of complex numbers ${\bf C}$.   The Gel'fand transform (which in the present case coincides with the Fourier--Laplace transform) maps $\mathscr{D} (R)$  to some algebra consisting of entire analytic functions on ${\bf C}$ according to the Payley--Wiener--Schwartz theorem.

The situation becomes especially transparent if $\Phi$ is a nuclear $F$-space (for example, $C^{\infty} (H)$ or $\mathscr{H} (H))$ or a complete nuclear $DF$-space (for example, $\mathscr{D} (H)$  or $\mathscr{A} (H))$. Then $\Phi$ is reflexive and $\mathscr{F} = \Phi'$  is a complete nuclear $DF$-space or nuclear $F$-space respectively.   The multiplication in $\Phi$ is automatically continuous and extends to a continuous map of the completed tensor product  $\Phi \hat{\otimes} \Phi$   into $\Phi$.   The dual map $\Phi' \to (\Phi \hat{\otimes} \Phi )' = \Phi' \hat{\otimes} \Phi'$ is called the comultiplication on $\Phi'$.   Analogously, the continuous linear map $\Phi \to \Phi \hat{\otimes} \Phi$ is called a comultiplication on $\Phi$ if the dual map induces a topological algebra structure on $\Phi'$.   One can define comultiplication in the standard way in the language of commutative diagrams, since passage to dual spaces and maps leads to "reversal of arrows" in diagrams and the tensor products $\Phi \hat{\otimes} \Phi$ and $\Phi' \hat{\otimes} \Phi'$  in the present case go into one another.   Thus, the diagram expressing the associativity of multiplication in $\Phi'$ goes into the diagram expressing the associativity of the comultiplication in $\Phi$, etc.   If there are given in $\Phi$ a continuous multiplication and comultiplication, then the pair $\mathfrak{A} = (\Phi , \Phi' )$  is a d.t.a. and conversely.   For example, if $H$ is an infinitely differentiable hypergroup, then the map $\varphi (x) \mapsto R^y \varphi (x)$ is a comultiplication $C^{\infty} (H) \to C^{\infty} (H) \hat{\otimes} C^{\infty} (H) = C^{\infty} (H \times H)$, and the multiplication of functions in $C^{\infty} (H)$ generates a comultiplication in $\mathscr{D} (H)$.

If the algebras $\Phi$  and $\Phi'$ have identities and the multiplications in $\Phi$ and $\Phi'$ are compatible in the sense that the comultiplication $\Phi \to \Phi \hat{\otimes} \Phi$ is a homomorphism of algebras (Hopf condition), then $\Phi'$ is called a topological Hopf algebra (t.H.a.).   In this case $\Phi'$ also has the structure of a t.H.a.   In the situation of Propositions 2 and 3 the Hopf condition means that $H$ is a semigroup.   In the presence of an involution $S$ this semigroup is a topological group (local compactness is not required).
\\

\underline{{\bf 3.6.}   Involution in Hypergroup Algebras.}   Let $H$ be a hypergroup with hypergroup algebra $\mathscr{F}$. If there is given in  $\mathscr{F}$ a map  $f \mapsto f^*$  such that $(f^*)^* = f$   and   $(f \ast g)^* = g^* \ast f^*$ for all $f, g \in H$,  then $\mathscr{F}$ is an {\it algebra with involution}; if $\mathscr{F}$ is a topological algebra, then the map $f \mapsto f^*$  is assumed to be continuous.  If $\mathscr{F}$ is an algebra with involution, then the right identity in $\mathscr{F}$ is also a two-sided identity, while the involution carries this element into itself.   It follows from this that  $\mathscr{F}$ coincides with its basic subalgebra, and the left and right generalized translations are completely equivalent.

There is a special interest in the case when the involution in $\mathscr{F}$ is induced by the involution  $x \mapsto x^*$  in the hypergroup $H$ (cf. above, point 2.7).   In this case in the space $\Phi$, where the g.t.o. act, there arises a map $\varphi (x) \mapsto \tilde{\varphi} (x) = \varphi (x^*)$   Combining this map with complex conjugation, one can define an involution in  $\mathscr{F}$ with the help of the equation
$$
\langle f^* , \varphi \rangle = \langle \overline{f, \tilde{\tilde {\varphi}}} \rangle . \eqno(3.5)
$$
(2.7) of point 2.7 means precisely that the involution defined by (3.5) is an antiautomorphism of the algebra $\mathscr{F}$,
i.e., $( f \ast g)^* = g^* \ast f^*$ for all $f, g \in \mathscr{F}$.

In particular, in group algebras there is a canonical involution, generated by passage to the inverse element in the group.   As indicated in point 2.7, reduced hypergroups, described in points 2.3-2.6, are hypergroups with involution.   Their hypergroup algebras are described in point 3.3 as subalgebras of group algebras.   It is easy to see that the involution in each of these hypergroup algebras is induced by the canonical involution in the corresponding group algebra.
\\

\underline{{\bf 3.7.}   Hypergroup Algebras, Formed of Functions.}   In order to define the generalized convolution of functions on the hypergroup $H$, it is necessary to fix a measure on $H$.   We return to the situation described in point 2.7:   the hypergroup $H$ is a locally compact space with a positive measure $m$.   With any function $a(x) \in L^1 (H, m)$ one can associate the bounded measure $a(x)m$ and define the convolution of functions as the convolution of the corresponding measures.   If $H$ is a hypergroup with involution, the measure $m$ is left-invariant, and the space $L^1 (H, m)$ is invariant with respect to real left generalized translations $L^x$, then the convolution $(a \ast  b)(x)$ of the functions $a(x)$ and $b(x)$ from $L^1 (H, m)$ is given by
$$
(a \ast b)(x) = \int a (y) \cdot L^{y^{*}} b(x)dm(y) . \eqno(3.6)
$$
It is clear that the corresponding analogous formula is valid for right-invariant measures and right generalized translation operators.

Let the measure $m$ be left-invariant and the involution in $H$ carry this measure into the equivalent right-invariant measure $m^* (x) = m(x^*)$.  We denote $\Delta (x)$ the modular function, i.e., the Radon--Nikodym derivative $dm(x) / dm^*(x)$.   The involution  $a(x) \mapsto a^*(x)$ induced in $L^1 (H, m)$ by the imbedding $L^1 (H, m) \to \mathscr{M}^b (H)$ described is given by
$$
a^* (x) = \Delta (x^*) \overline{a(x^*)} . \eqno(3.7)
$$
To study hypergroup algebras of the type of $L^1 (H, m)$ one can use the powerful apparatus of the theory of Banach algebras with involutions (cf. Introduction).

By an analogous scheme one also constructs other hypergroup algebras consisting of functions (for example, finite ones).   One can define a hypergroup  $C^*$-algebra as the enveloping $C^*$-algebra of the algebra $L^1 (H, m)$.   In [17], in connection with the proof of the Plancherel theorem for hypergroups, hypergroup left Hilbert algebras were considered (cf., e.g., [82] on left Hilbert algebras).   The scalar product in this Hilbert algebra is induced by the imbedding in $L^2 (H, m)$, and the multiplication (convolution) and involution are defined in terms of the adjoint g.t.o.   For hypergroups with involution and real g.t.o. these definitions reduce to (3.6) and (3.7).

\section*{\small 4.     Infinitesimal   Hypergroup   Algebras   and   Generalized   Lie   Theory}

\underline{{\bf 4.1.}   Laplace Transform.}   Let $\Phi$ be a space of functions (germs, formal power series) of the variables $x_1, \ldots , x_n$, containing exponentials of the form $\exp (\zeta_1 x_1 + \ldots + \zeta_n x_n)$.   By the {\it Laplace transform} of the functional $\int \colon~ \varphi \mapsto \langle f, \varphi \rangle$ defined on $\Phi$ is meant the function $f(\zeta_1 , \ldots , \zeta_n ) = \langle f, e^{\zeta_1 x_1 , \ldots , \zeta_n x_n} \rangle$.

If $\Phi$ coincides with the space $\mathscr{H}_0$ of germs of analytic functions at the point $x_1 = \ldots = x_n =0$, endowed with the standard toplogy\symbolfootnote[1]{The sequence $\varphi_n$ converges in this topology if it converges uniformly on some complex neighborhood of the point zero, where all the $\varphi_n$ are defined in this neighborhood.}, then the dual space $\mathscr{F}_0$  of $\mathscr{H}_0$  can be identified with the help of the Laplace transform with the space of entire functions, which grow slower than $\exp [\varepsilon (| \zeta_1 | + \ldots | \zeta_n | )]$ for arbitrarily small $\varepsilon > 0$.   For all $f \in \mathscr{F}_0 , \varphi \in \mathscr{H}_0$ one has
$$
\langle f, \varphi \rangle = \sum_{s_1 , \ldots , s_n} \frac{f_{s_1 , \ldots , s_n} \varphi_{s_1 , \ldots , s_n}}{s_1 ! s_2 ! \ldots s_n ! }, \eqno(4.1)
$$
where $\varphi = \sum_{s_1 , \ldots , s_n} \varphi_{s_1 , \ldots , s_n} \frac{x_1^{s_1} \ldots x_n^{s_n}}{s_1 ! \ldots s_n !} , f(\zeta_1 , \ldots , \zeta_n )  = \sum_{s_1 , \ldots , s_n} f_{s_1 , \ldots , s_n} \frac{\zeta_1^{s_1} \ldots \zeta_n^{s_n}}{s_1 ! \ldots s_n !}$ are the expansion in a Taylor series of the germ of $\varphi$ and the Laplace transform of the functional $f$.

If $\Phi$ coincides with the space of germs at zero of infinitely differentiable functions of the real variables $x_1 ,\ldots x_n$ (endowed with the standard topology of uniform convergence along with derivatives on some neighborhood of zero), then the dual space consists of all generalized functions with support at zero. The Laplace transform identifies this space with the space ${\bf C} [\zeta_1 , \ldots , \zeta_n ]$ of all polynomials of $\zeta_1 , \ldots , \zeta_n$; (4.1) remains valid.

Finally, let $\Phi$ coincide with the space ${\bf C}[[x_1 , \ldots , x_n ]]$ of formal power series of $x_1 ,\ldots , x_n$.   If we endow $\Phi$ with the strongest locally convex topology, then $\Phi'$ consists of all linear functionals on $\Phi$.   The Laplace transform identifies $\Phi'$ with the space of polynomials ${\bf C} [ \zeta_1 , \ldots , \zeta_n ]$, and the duality between ${\bf C} [[x_1 , \ldots , x_n ]$ and ${\bf C}[\zeta_1 , \ldots , \zeta_n ]$ is given by (4.1).
\\

\underline{{\bf 4.2.}   Algebras Associated with  Spaces of Germs.}   We consider a double topological algebra  $\mathfrak{A} =(\mathscr{H}_0 , \mathscr{F}_0)$ where $\mathscr{H}_0$ is the algebra of germs of analytic functions at zero (cf. above) with respect to multiplication, $\mathscr{F}_0$ is the dual space of $\mathscr{H}_0$. If the multiplicative functional $\varphi \mapsto \varphi (0)$  on $\mathscr{H}_0$   is a right identity in the algebra $\mathscr{F}_0$, then by analogy with the results of points 3.4 and 3.5 it is natural to consider $\mathscr{F}_0$ as the hypergroup algebra of a "local hypergroup".   In the one-dimensional case such algebras were studied in [42].   Algebras of the type indicated arise naturally as subalgebras of hypergroup algebras of analytic functionals.   Associative ultra-envelopes of Lie algebras in the sense of P. K. Rashevskii [77] correspond to the group case.
\\

\underline{{\bf 4.3.}   Formal Hypergroups and Generalizations of the Poincare --}\\ \underline{Birkhoff -- Witt Theorem.} By virtue of the duality between ${\bf C}[[x_1 , \ldots , x_n]]$ and ${\bf C}[\zeta_1 , \ldots , \zeta_n ]$ the associative comultiplication in ${\bf C}[[x_1 , \ldots , x_n ]]$ induces an associative multiplication in ${\bf C}[\zeta_1 , \ldots , \zeta_n ]$ (which in general does not coincide with the usual multiplication of polynomials), so that a dual algebra arises.   If it has a right coidentity, which is a multiplicative functional (i.e., the polynomial $\varepsilon \equiv 1$ is a right identity with respect to the induced multiplication), then by analogy with the results of points 3.4 and 3.5, it is natural to consider that the comultiplication in ${\bf C}[[x_1 , \ldots , x_n]]$ defines the structure of a formal hypergroup, and the dual algebra is a "hypergroup" one.   It is also natural to require that the hypergroup multiplication be compatible with the natural filtration in ${\bf C}[\zeta_1 , \ldots , \zeta_n ]$ and that the identity in ${\bf C}[[x_1 , \ldots , x_n ]]$ be a multiplicative functional on ${\bf C}[\zeta_1 , \ldots , \zeta_n ]$.  Obvious generalizations are connected with replacement of the field $C$ by other fields and rings, passage to algebraic extensions of algebras of power series, the introduction of anticommuting variables, etc.   Projectors generated by the right identity in the hypergroup algebra single out the basic subalgebra in this algebra and the basic subspace in ${\bf C}[[x_1 , \ldots , x_n ]]$ by analogy with the constructions of points 2.3 and 3.2.   The dual algebra which arises defines a reduced formal hypergroup.   Formal hypergroups and generalizations of them from various points of view arise in many questions of topology, analysis, and mathematical physics, cf., e.g., [9, 11-13, 19, 28-30, 33, 34, 36, 37, 52-56, 63, 65-67, 79-81, 88, 108, 162].

In the group case the scheme indicated leads to universal enveloping algebras of Lie algebras.   Let $\mathfrak{g}$ be a Lie algebra with basis $X_1 , \ldots , X_n$.   The enveloping algebra   $U (\mathfrak{g})$   is the algebra of polynomials with complex coefficients of the generators $X_1 , \ldots , X_n$, connected by the relations $X_i X_j = X_j X_i + \sum_{k=1}^n C^k_{ij} X_k$, where $C^k_{ij}$ are the structural constants of the Lie algebra $\mathfrak{g}$.    The Poincare--Birkhoff--Witt theorem (for short PBW) lets one construct a one-to-one linear correspondence $\Lambda \colon U(\mathfrak{g} ) \to {\bf C} [\zeta_1 , \ldots , \zeta_n ]$  compatible with the filtrations in these algebras.   As a result the structure of associative algebra can be carried over from  $U(\mathfrak{g} )$ to ${\bf C}[\zeta_1 , \ldots , \zeta_n ]$, and by duality there arises in ${\bf C} [[x_1, \ldots , x_n]]$ a (group) multiplication.   The construction of the map $\Lambda$ is connected with the choice of local coordinates $x_1 , \ldots , x_n$ in the Lie group $G$ with the given Lie algebra $\mathfrak{g}$.   By L. Schwartz theorem $U(\mathfrak{g})$ can be identified with the subalgebra of the group algebra $\mathscr{D} (G)$   consisting of all generalized functions with support at the identity of the group.   The Laplace transform in canonical coordinates of the second kind connected with the basis $X_1, \ldots , X_n$ in $\mathfrak{g}$ carries the element $\sum f_{s_1 , \ldots , s_n} X_1^{s_1} \ldots X_n^{s_n} \in U(\mathfrak{g})$ (more precisely, its image in $\mathscr{D}(G)$) into the polynomial $\sum f_{s_1 , \ldots , s_n} \zeta_1^{s_1} \ldots \zeta_n^{s_n}$; the Laplace transform in canonical coordinates of the first kind carries the element $\sum f_{s_1 , \ldots , s_n} ((X_1^{s_1} \ldots X_n^{s_n}))$,  where $((X_1^{s_1} \ldots X_n^{s_n}))$   is the result of symmetrization of the monomial $X_1^{s_1} \ldots X_n^{s_n}$, into this polynomial (cf. [59, Chap. III]).
 
For many associative algebras with identity generated by a finite collection of generators $X_1 ,\dots , X_m$ (and classes of such algebras), it is proved (or assumed) that the analogs of the PBW theorem are valid (cf. [19, 28-30, 33, 34, 36, 37, 65, 66, 80, 81]); one can also give analogs of the Jacobi identity, guaranteeing the validity of a PBW type theorem in this case or that.   For algebras generated by quadratic commutation relations, possible generalizations of the PBW theorem are analyzed in [19].   Algebras for which analogs of the PBW theorem are valid are natural candidates for the role of formal and infinitesimal hypergroup algebras.
\\

\underline{{\bf 4.4.}   Infinitesimal Hypergroup Algebras.}   Let $H$ be an infinitely differentiable hypergroup with hypergroup algebra $\mathscr{D}(H)$   consisting of generalized functions with compact support (cf. above, point 3.4).   By the infinitesimal hypergroup algebra of the hypergroup $H$, we mean the subalgebra $\mathscr{D}_e (H)$   generated by generalized functions with support at the identity $e \in H$.   If the support of any functional from $\mathscr{D}_e (H)$  coincides with $e$, then the Laplace transform in any local coordinates maps $\mathscr{D}_e (H)$  one-to-one and linearly onto the space of polynomials, and we return to the situation considered in point 4.3.

The correspondence $f \mapsto R(f)$  is a representation of the algebra $\mathscr{D}_e (H)$ on the space $C^{\infty} (H)$ and maps this algebra onto the algebra of infinitesimal right translations $\mathscr{D}_e^R (H)$;   analogously, the antirepresentation $f \mapsto L(f)$   maps $\mathscr{D}_e (H)$ onto the algebra of infinitesimal left translations $\mathscr{D}_e^L (H)$.   The algebras $\mathscr{D}_e (H)$ and $\mathscr{D}_e^L (H)$  are isomorphic.   If $\delta_c$ is a two-sided identity in $\mathscr{D}_e (H)$ then the algebras $\mathscr{D}_e (H)$ and  $\mathscr{D}_e^R (H)$ are isomorphic.   In general the algebra $\mathscr{D}_e^R (H)$  is isomorphic to the intersection of $\mathscr{D}_e (H)$  with the basic subalgebra of $\mathscr{D}_e (H)$; hence $\mathscr{D}_e^R (H)$ should be considered as the infinitesimal hypergroup algebra of the reduced hypergroup (cf. above points 2.3, 3.2, and 3.4), although the reduced hypergroup is not necessarily a smooth manifold.

The construction of infinitesimal hypergroup algebras and the reconstruction of "global" hypergroups from given "infinitesimal" objects is a natural subject of Lie theory.   The construction of such a theory is far from to be completed, although important results in this direction were already found in the first papers devoted to g.t.o.
\\

\underline{{\bf 4.5.}   Generators and Lie Theorems for Hypergroups.}   Let $H$ be an infinitely differentiable hypergroup, so that $u(x, y) = R^x \varphi (y)$ is a smooth function on $H \times H$ for all $\varphi \in C^{\infty} (H)$.    We denote by $x_1 , \ldots , x_n$ local coordinates of a point $x \in H$ and we choose a coordinate system such that the neutral element has zero coordinates. By generators of right translation of order $n$ we mean linear operators of the form
$$
\left. R_{s_1 , \ldots , s_n ;h}~ \colon~ \varphi (h) \mapsto \frac{\partial^s u (x,h)}{\partial x^{s_1}_1 \ldots \partial x^{s_n}_n} \right|_{x=0} , \eqno(4.2)
$$
where $u(x, h) =R^x \varphi (h), s = s_1 + \ldots + s_n$.   Analogously, generators of left translation have the form:
$$
\left. L_{s_1 , \ldots , s_n ;h}~ \colon~ \varphi (h) \mapsto \frac{\partial^s u (x,h)}{\partial x^{s_1}_1 \ldots \partial x^{s_n}_n} \right|_{x=0} . \eqno(4.3)
$$
In exactly the same way one can define generators for complex local coordinates on a holomorphic hypergroup for $\varphi \in \mathscr{H} (H)$.

It follows from the associativity axiom that any generator of left translation commutes with any generator of right translation (just as with the operators $R^x$).  Differentiating the associativity condition  $R^x L^y \varphi (h) = L^y R^x \varphi (h)$ with respect to the variables $h_1 , \ldots , h_n$  the corresponding number of times, and then letting $h_1 = \ldots = h_n =0$, we get the system of equations
$$
L_{s_1 , \ldots , s_n ;x} (u) = R_{s_1 , \ldots , s_n ;y} (u) , \eqno(4.4)
$$
where $u(x,  y) = R^y \varphi (y)$.   One should consider the system (4.4) as a generalization of the first Lie theorem [52, 54-56].   For Delsarte hypergroups this result is found in [108].   For the hypergroup described in point 2.6, where $G$ is a semisimple linear group   and $K$ is its maximal compact subgroup, a complete collection of generators is actually calculated in [23, 10].

It is not necessary to use all the equations of the system (4.4) to determine $u(x, y)$.   For example  for a translation on a Lie group, the first order generators already determine the function u uniquely    (i.e., the group multiplication).   In general some generators of lower orders may degenerate (for example, to multiplication by a constant), so that the corresponding equations of the system (4.4) do not contain useful information.  Hence an important problem arises: to select the minimal number of equations from the system (4.4), which uniquely determine the g.t.o. Here the degenerated generators enlarge the number of initial conditions. If a finite system of the form (4.4) under certain initial conditions including the condition  $u_{1 x=0} = \varphi (y)$  uniquely determines the solution $u(x, y) = R^x \varphi (y)$, where the operators on the left sides of the system commute with all the operators on the right sides, then $R^x$ are g.t.o. (i.e., $H$ is a hypergroup). This assertion is an analog of the first Lie inverse theorem [53-56].   For a certain class of g.t.o. B.M. Levitan proved analogs of the second and third (direct and inverse) Lie theorems. In particular, in the space of smooth functions of n variables, g.t.o.  are constructed for which the generators of right (left) translation generate any given n-dimensional Lie algebra [52 56]. In [25, 126] a complete explicit description of these generators is found in the form of second order integrodifferential operators.   With the help of an analogous technique generators of any order acting on the space of entire analytic functions of $n$ variables and generating any n-dimensional real Lie algebra are constructed in [24]; from these generators one can reconstruct the hypergroup.   One can construct hypergroups starting from commutation relations of more general type [53, 65-67, 36, 37, 28 30].

We note that already in Delsarte's first paper [105] g.t.o. on the line were constructed, starting from an explicitly given second order generator, with the help of a series analogously to the way in wich an ordinary translation decomposes into a series in powers of the differentiation operator (Taylor series).  B. M. Levitan and A. Ya.  Povzner studied commutative g.t.o. on the line with second order generator of Sturm --  Liouville type (cf.   [51   55   56]); a complete classification of one-dimensional infinitely differentiable hypergroups with Sturm -- Liouville type generators (including noncommutative ones) is given in [27]. A vast literature is devoted to g.t.o. on the line generated by a differential operator, cf., e.g., [7, 51, 55, 56, 83, 90, 97, 98, 109, 134, 146].

It is easy to verify that the generators (4.2) and (4.3) can be represented in the form $R(f)$ and $L(f)$ respectively, where the element $f \in \mathscr{D} (H)$ has the form
$$
\left. f \  \colon \  \varphi \mapsto \langle f, \varphi \rangle =  \frac{\partial^s \varphi (x)}{\partial x^{s_1}_1 \ldots \partial x^{s_n}_n} \right|_{x_1 = x_2 =\ldots = x_n = 0} . \eqno(4.5)
$$
Since any generalized function with support at the point $e \in H$ is a finite linear combination of the element $\delta_e$ and generalized functions of the form (4.5), the following assertion is valid.
\\

\underline{Proposition 4.}   Generalized functions of the form (4.5) and the delta-function $\delta_e$ generate the infinitesimal hypergroup algebra $\mathscr{D}_e (H)$, and the generators of left (right) translation and the operator $L^e$ (respectively $R^e =\mathfrak{D}$ generate the algebra of infinitesimal left (right) translations $\mathscr{D}_e^L (H)$  (respectively   $\mathscr{D}_e^R (H)$).
\\

\underline{{\bf 4.6.}  Example.}   We return to the example described in point 2.4.  Differentiating (2.3), it is easy to calculate the generators $R_n$ and $L_n$ of left and right translation of the n-th order:   $R_1 = R_3 = \ldots = R_{2k+1} = \ldots = 0; R_2 = \frac{d^2}{dt^2} , R_{2k} = (R_2)^k; L_n = L_0  \frac{d^n}{dt^n}$ where $L_0= L^e  \varphi (t) \mapsto [\varphi (t) + \varphi (-t)] /2$.

The algebra of infinitesimal right translations $\mathscr{D}^R_0$,    coincides with the algebra of polynomials in the operator $R_2 = d^2 / dt^2$ and is isomorphic to the subalgebra, consisting of even generalized functions with support at zero, of the infinitesimal hypergroup algebra $\mathscr{D}_0$.     For the elements of this subalgebra, the generalized convolution coincides with the ordinary one; the Laplace transform $f \mapsto f(\zeta )$  carries the indicated subalgebra into the algebra of polynomials of $\zeta^2$ with ordinary multiplication.

The algebra of infinitesimal left translations $\mathscr{D}^L_0$ and the algebra $\mathscr{D}_0$, antiisomorphic to it, have a more complex nature.   The operators $L_n = L_0 \cdot (d^n / dt^n)$ form a basis in $\mathscr{D}^L_0$ while  $L_{2k+1} L_n = 0$  and $L_{2k} L_n = L_{2k+n}$ for all $n, k = 0, 1, 2, 3, \ldots$; it follows from this that $L_{2k} = (L_2 )^k$ and $L_{2k+1}  = (L_2 )^k L_1$.   Thus, the elements $L_0$ (left identity), $L_1$ and $L_2$ generate the algebra $\mathscr{D}^L_0$  (here $L_1 L_0 = L^2_1 = L_1 L_2 =0, L_2 L_0 =L_2 , L_2 L_1 = L_3 , L_2^2 =L_4$), and the analog of the Poincare--Birkhoff--Witt theorem obviously holds.   The structure of the algebra $\mathscr{D}_0$  is easily determined by virtue of the antiisomorphism between $\mathscr{D}_0$ and $\mathscr{D}_0^L$.

The system of equations (4.4) in the present case reduces to the wave equation $\partial^2 u / \partial x^2 = \partial^2 u / \partial t^2$ and the function $u(x, t) = R^x  \varphi (t) = [\varphi (t+x)+ \varphi (t-x)] /2$ is a solution of this equation with initial conditions $u(0, t) = \varphi (t)$ and $u'_x (0, t) = 0$.  On the reduced hypergroup, i.e., on the half-line $0 \leq t < \infty$, one can define generalized translation as a solution of the wave equation for $t \geq 0$ with the same initial conditions.   If the generator $d^2 / dt^2$ is replaced by a Sturm--Liouville operator with potential $v(t)$, then the corresponding generalized translation on the half-line $0 \leq t < \infty$ can be defined as a solution of the equation
$$
\frac{\partial^2 u}{\partial x^2} - \frac{\partial^2 u}{\partial t^2} + v(x) - v(t) = 0
$$
with the initial conditions indicated above.
\\

\underline{{\bf 4.7.}   Example.}   We consider the group $G = SU(2)$, and in the Lie algebra of this group we choose a standard basis $X, Y, Z$, such that $[X, Y] = Z,~ [Y, Z] =X,~ [Z, X] = Y$.   Let $a$ and $b$ be automorphisms of the group $G$, carrying the basis $(X, Y, Z)$ respectively into $(-X, Y, -Z)$ and $(X, -Y, -Z)$.   The equation
$$
R^x \varphi (t) = \frac{1}{2} [\varphi (t \cdot x) +\varphi (b(t) \cdot x) + \varphi (t \cdot a(x)) - \varphi (b(t) \cdot a(x))]
$$
defines a g.t.o. in $C^{\infty} (G)$, and the first order generators of right translation $X, Y$, and $Z$ satisfy the relations $XY + YX =Z,~ YZ + ZY =-X,~ ZX + XZ =-Y$ (cf. [28]).

It is easty to verify that the algebra with identity generated by these generators is isomorphic to the infinitesimal hypergroup algebra and the analog of the PBW theorem is valid.

\section*{\small 5.     Special   Classes   of  Hypergroups}

\underline{{\bf 5.1.}   Generalization of Hopf's Condition.}   The contents of Section 3 (in particular, Theorem 2 and its analogs) show that hypergroups form an arbitrarily large class, that to get meaningful results it is reasonable to include additional restrictions and to consider special classes of hypergroups.   For example, one can single out classes of hypergroups in terms of generalized Lie theory (cf. above).   It is also natural to require that the comultiplication $\varphi (x) \mapsto \psi (x,y) = R^y \varphi (x)$ be compatible to some degree or other with ordinary multiplication of functions.   A relatively weak requirement is that the space $\Phi$, where the g.t.o. act, is a (commutative) algebra with identity, while the identity in $\Phi$ is a multiplicative functional on the hypergroup algebra $\mathscr{F} = \Phi'$ and the identity in  $\mathscr{F}$ is a multiplicative functional on $\Phi$.   The stronger Hopf condition says that the hypergroup reduces to a semigroup (cf. above, point 3.5).   However even for hypergroups which do not reduce to semigroups, one is sometimes able to give $\Phi$ the structure of a (noncommutative) algebra in such a way that the comultiplication is a homomorphism $\Phi \to \Phi \widehat{\otimes} \Phi$.    Infinitesimal and formal objects of this kind, called generalized Lie algebras and groups, are introduced in [28, 30]; here the concept of tensor product is generalized in the spirit of MacLane [147] with the help of a symmetry operator (cf. also [63, 67]).   Generalized Lie groups and algebras are closely connected with Lie supergroups and superalgebras, and also with the Yang--Baxter equation which is popular in theoretical physics.   The hypergroup described in point 4.7  is an example of a (global) generalized Lie group; one can construct a new associative operation in $C^{\infty} (G)$ by combining multiplication of functions with the action of the automorphisms $a$ and $b$.

Some other conditions on the comultiplication, replacing the Hopf condition, are considered below.
\\

\underline{{\bf 5.2.}   Multivalued Groups.}   Let $h$ be a hypergroup and for any elements $x, y \in H$ let the generalized convolution of the delta-functions $\delta_x$ and $\delta_y$ have the form
$$
\delta_x \ast \delta_y = \sum^n_{i=1} a_i \delta_{x_i} ,
$$
where $a_i$ are numerical coefficients depending on $x$ and $y$, and $x_i$ are elements of $H$, depending on $x$ and $y$.  In this case, following [29] (where differentiable and holomorphic hypergroups are considered), we call $H$ a metamultivalued group; if $a_1 = a_2 = \ldots = a_n = 1/n$, then $H$ is called an n-valued group (the connection with the definition of abstract hypergroup in point 1.2 is obvious).  An example of a two-valued hypergroup is described in point 2.4, and one of a metamultivalued group is described in point 4.7.   In [29] regular methods of construction of hypergroups of this type are given.   One-dimensional formal multivalued groups are introduced by V. M. Bukhshtaber in connection with topological applications, cf. [12] for a surveying account; there is further development along these lines in [13, 88], for example.
\\

\underline{{\bf 5.3.}   Hypergroups with Positivity Conditions.}   An enormous literature is devoted to hypergroups with involution, for which the Hopf condition is replaced by the requirement that the generalized convolution have certain positivity properties.   Of this class are p-hypergroups and hypercomplex systems (h.s.) with continuous bases (cf. Introduction).   Similar objects (of more general nature) such as dual algebras with involution and coinvolution, identity and coidentity, and with additional positivity conditions  were considered by A. M. Vershik [18] in connection with the geometric theory of $C^*$-algebras and duality in representation theory; cf. also [39]. At the present time the use of p-hypergroups (usually called simply hypergroups\symbolfootnote[2]{ln [136], which contains the "canonical" axiomatics of $p$-hypergroups, the term "convo" (from the word convolution) is used.}) is of the greatest popularity.

Let $H$ be a locally compact hypergroup with involution (cf. above, points 2.7, 3.2, and 3.6), where the space $C_0 (H)$ of continuous functions which vanish at infinity  is invariant with respect to the g.t.o. and the involution, so that not only the space $\mathscr{M} (H)$   of measures with compact support   but also the space $\mathscr{M}^b (H)$ of bounded measures is endowed with the structure of a hypergroup algebra with involution.   We denote by $\mathscr{M}^b_{+} (H)$ the set of all nonnegative measures from $\mathscr{M}^b (H)$ endowed with the weakest topology for which the maps $f \mapsto f (H) = \int_H df$  and $f \mapsto \langle f, \varphi \rangle = \int \varphi df$  are continuous for all nonnegative continuous functions $\varphi$ with compact  support in $H$.

The hypergroup $H$ is a $p$-hypergroup if the following conditions hold:   1) the generalized convolution of nonnegative measures is a nonnegative measure, and the convolution induces a continuous map $\mathscr{M}^b_+ (H) \times \mathscr{M}^b_+ (H) \to \mathscr{M}^b_+ (H)$;   2) for any $x, y \in H$ the convolution of the Dirac measure (delta-functions) $\delta_x$ and $\delta_y$ is a probability measure, so that $\delta_x \ast \delta_y (H) = \int d(\delta_x \ast \delta_y )$; 3) the support supp$(\delta_x \ast  \delta_y )$ of the measure $\delta_x \ast \delta_y$ contains the neutral element (identity) of the hypergroup if and only if $x = y^*$; 4) the map $(x, y) \mapsto$ supp$(\delta_x \ast \delta_y )$   is a continuous map from $H \times H$ to the space $\mathfrak {R} (H)$  consisting of compact subsets of $H$ and endowed with the Michael topology [154]; this topology is generated by the subbase ${K \in \mathfrak{R} (H)   \colon K \cap U \ne \varnothing , K \subseteq V}$, where $U$ and $V$ run through the collection of open subsets of $H$.   The present definition is a slight modification of the definition formulated in [136]; for other versions of the axiomatics of p-hypergroups, cf. [112, 113, 128, 169, 170, 160]; in [128] even non-Hausdorff hypergroups, which arise as dual objects to certain groups, are considered.  The hypergroups described in points 1.2, 1.3, 2.3-2.6 (cf. also points 2.7, 3.3, and 3.6)   are examples of p-hypergroups.  Other interesting examples are described, for example, in [94, 115, 127-132, 142-144, 160, 161, 169, 182, 183].  These hypergroups can also be treated as h.s. with locally compact basis.

Objects of this type are similar to locally compact groups in their properties to a large degree.   The theory of commutative h.s. and $p$-hypergroups is especially well developed.   The basic results of harmonic analysis, including Pontryagin duality and refined questions of spectral analysis and synthesis, generalize to this case (cf. Introduction).  On commutative h.s. and p-hypergroups one can construct the nuclear space of "basic" finite functions  with the invariance property [4].   The existence is proved of a positive invariant measure on commutative p-hypergroups [170]; cf. [8, 35] for the corresponding result for h.s. (existence of a multiplicative measure).   The existence of an invariant measure is proved for compact and discrete $p$-hypergroups [136, 169]; the question is open for the general case \footnote{However the uniqueness of a left-invariant (right-invariant) measure is proved up to a constant factor.}.   If an invariant measure m exists, then one can construct the hypergroup algebra $L^1 (H, m)$ (cf. above, point 3.7) and use the methods of the theory of Banach algebras with involution.

By no means all interesting hypergroups with involution are $p$-hyper-groups or h.s. with locally compact bases.   A standard example is the hypergroup generated by a Sturm -- Liouville operator on the half-line with rapidly growing potential (cf. above, points 1.5 and 4.6).   Even in the compact case passage to the dual hypergroup leaves the class of $p$-hypergroups (h.s.) (for example, cf.  [3, 7, 112, 136, 169, 182]) due to the loss of the positivity condition.   Nevertheless the class of $p$-hypergroups (not necessarily commutative) can be immersed in a larger class (consisting of objects of the type of Kac algebras), for which a satisfactory duality theory is constructed, cf. [138].
\\

\underline{{\bf 5.4.}   Gelfand Pairs.}   Let  $H = K \setminus G / K$ be the hypergroup of double cosets of the locally compact group $G$ with respect to its compact subgroup $K$ (cf. above, points 2.6, 2.7, 3.3, and 3.6).   If the hypergroup $H$ is commutative, then the pair $(G, K)$ is called a {\it Gelfand pair.} These objects were studied intensively (in particular, in the framework of the theory of p-hypergroups and h.s.), among others in connection with the theory of spherical functions and the theories of probability measures on groups and hypergroups (for example, in connection with the generalization of the central limit theorem to the case of convolution of probability measures), cf. [23, 7, 10, 94, 96, 110, 111, 125, 131, 132, 136, 137, 145, 151, 167-169, 175, 196].   The following result goes back to Gel'fand's work [23], cf. also [131].
\\

\underline{GELFAND'S THEOREM.}   Let the group $G$ be unimodular, and the compact subgroup $K$ coincide with the set of fixed points of the involution automorphism $\sigma$.   If any element $x \in G$ admits a decomposition $x = k \cdot y$, where $k \in K, \sigma (y) = y^{-1}$. then $(G, K)$ is a Gelfand pair.
\\

It is easy to see that to prove this theorem it suffices to establish the commutativity of the hypergroup algebra $L^1 (H)$.   We identify $L^1 (H)$ with a subalgebra of the group algebra $L^1 (G)$ (cf. points 3.6 and 3.7); this subalgebra consists of all K-biinvariant functions in $L^1 (G)$.   Let $f(x)$ be such a function.   We set $f^{\sigma} (x) = f(\sigma (x)), \tilde{f}(x) = f(x^{-1})$. Then $f^{\sigma} (x) = f(\sigma (x)) = f(\sigma (k \cdot y)) = f(\sigma (k) \cdot \sigma (y)) = f(k \cdot y^{-1}) = f(y^{-1} \cdot k^{-1}) = f (x^{-1}) = \tilde{f} (x)$. The map $f \mapsto f^{\sigma}$  is an automorphism of the algebra $L^1 (H)$ and preserves the order of the factors in the product of elements of $L^1 (H)$; the map $f \mapsto  \tilde{f}$  is an antiautomorphism and changes the places of the factors in the product of two elements of $L^1 (H)$.  Since these maps coincide, the algebra $L^1 (H)$ is commutative, which is what was needed.

The hypothesis of Gelfand's theorem holds, if $G$ is a semisimple linear group or group of motions, and $K$ is a maximal subgroup of it, and also if $(G, K)$ is a Riemannian symmetric pair in the sense of [84, Chap. VI, Sec. 1].

\section*{\small 6.    Representations   of Hypergroups  and   Harmonic   Analysis}

\underline{{\bf 6.1.}  Representations of Hypergroups and Hypergroup Algebras.}   The theory of representations of hypergroups is to a large extent analogous to the theory of representations of groups.   This analogy extends quite far.   For example, with the help of such concepts as generator and infinitesimal hypergroup algebra, one can study representations of differentiable and holomorphic hypergroups by an infinitesimal method just as is done for Lie groups; cf. [108, 55, 56, 24], in which finite-dimensional representations of Delsarte hypergroups (cf. above, point 2.3) and hypergroups generated by generators of the second order, generating the Lie algebra (cf. point 4.5) are considered.   It is convenient to treat representations of hypergroups as representations of hypergroup algebras.  Just as in the case of groups, to different types of hypergroup algebras correspond different versions of the theory of representations and harmonic analysis (cf. [32, 59, 61]).

Let $V$ be an arbitrary locally convex space, $S(V)$ be the algebra of all weakly continuous linear operators  on $V$, endowed with the weak operator topology (cf. above, point 3.1), $\mathscr{F}$ be an arbitrary topological algebra  (cf. point 3.2).   By a continuous representation of the algebra $\mathscr{F}$ on the space $V$ we shall mean a continuous map $\mathscr{F} \to S(V)$  which is a homomorphism of algebras; if $\mathscr{F}$ has a two-sided identity, we shall assume (without special mention)   that the identity goes into the identity operator.

Let $H$ be a hypergroup, where the corresponding g.t.o. $R^x$ and $L^x$ act on the locally convex space $\Phi$.  By a representation of the hypergroup $H$ on the space $V$ we mean a map $x \mapsto \pi_x$ of the space $H$ into $S(V)$, such that the following conditions hold:

1)   for any elements $v \in V, v' \in V'$ the function $\varphi_{v, v'} (x) = \langle v', \pi_x v \rangle$, called a {\it matrix element of the representation}, is contained in $\Phi$ and  $R^y \varphi_{v,v'} (x) = \langle v', \pi_x  \pi_y v \rangle$ for all $x, y \in H$;

2)   for any element $f$ of the dual space $\Phi'$ to $\Phi$  there exists an operator $\pi (f) = \int \pi_x df (x) \in S(V)$.

It follows from the definition of the integral $\int \pi_x df(x)$ (cf. point 3.1), Eq. (3.1), and the fact that the matrix elements of the representation $\pi$ are contained in $\Phi$,  that the map $f \mapsto \pi (f) = \int \pi_x df(x)$  of the space $\Phi'$ into $S(V)$ is continuous in the weak topology $\sigma (\Phi' , \Phi)$ and consequently  in the Mackey topology $\tau (\Phi' , \Phi)$ and in any topology compatible with the duality between $\Phi'$ and $\Phi$.   If the generalized convolution turns $\Phi'$ into a topological hypergroup algebra, then the map $f \mapsto \pi (f)$  is a continuous representation of the algebra $\Phi'$.   Condition 2) holds automatically, if $V$ is a complete or quasicomplete barreled space\symbolfootnote[1]{We recall that in this case any operator from $S(V)$ is continuous.} (for example, a reflexive or complete metrizable space), and $\Phi$ coincides with the space $C(H)$ of continuous functions ($C^{\infty} (H)$ or $\mathscr{H} (H)$ of infinitely differentiable or holomorphic functions) on the locally compact (respectively infinitely differentiable  or holomorphic) hypergroup $H$, cf. [61].   We call a representation $\pi$ of the hypergroup $H$ {\it nondegenerate}, if the  operator $\pi_e$ (where $e$ is the neutral element in $H$) is the identity operator.   If $\Phi$ coincides with its basic subspace (cf. point 2.1), then we shall assume representations of the hypergroup $H$ to be nondegenerate without any special mention.

Let $H$ be a locally compact hypergroup.   Any representation $x \mapsto \pi_x$ of it extends to a continuous representation $f \mapsto \pi (f)$ of the hypergroup algebra $\mathscr{M} (H)$   (cf. above).  Conversely, if we have a continuous representation $f \mapsto \pi (f)$ of the algebra $\mathscr{M} (H)$,    then the correspondence $x \mapsto \pi_x = \pi (\delta_x)$ is (as is easy to verify, cf. [61]) a representation of the hypergroup $H$.   It is easy to derive the following theorem from this.
\\

\underline{THEOREM 3.} The representations of a locally compact hypergroup $H$ are in one-to-one correspondence with continuous representations of the topological hypergroup algebra $\mathscr{M} (H)$ . If a representation of the hypergroup $H$ is nondegenerate, then the corresponding representation of $\mathscr{M} (H)$ is the composition of a continuous representation of the basic subalgebra\symbolfootnote[2]{Cf. point 3.2 for the definition of the basic subalgebra.} $\widetilde{\mathscr{M}} (H)$ and the canonical projection $\mathscr{M} (H) \to \widetilde{\mathscr{M}} (H)$. Thus a one-to-one correspondence is established between nondegenerate representations of the hypergroup $H$ and continuous representations of the algebra   $\widetilde{\mathscr{M}} (H)$.
\\

Analogous results are valid for infinitely differentiable and holomorphic hypergroups.   To representations with infinitely differentiable (holomorphic) matrix elements correspond continuous representations of hypergroup algebras of generalized functions with compact support (algebras of analytic functionals), while to non-degenerate representations of hypergroups correspond representations of basic subalgebras, cf. [61].

The analog of the theory of unitary representations of groups is the theory of symmetric (Hermitian) representations of hypergroup algebras with involution on Hilbert spaces (the representation $\pi$ is {\it symmetric} if the operator $\pi (f^*)$  is adjoint to the operator $\pi (f)$  where $f \mapsto f^*$  is the involution in the hypergroup algebra). Representations of this type were studied (primarily for the commutative case) by B. M. Levitan, Yu. M. Berezanskii, S. G. Krein, and other authors, cf., e.g., [3, 5, 7, 8, 14, 17, 47-51, 93, 96, 112, 129, 131, 132, 135, 136, 142, 144, 149, 150, 169, 180].   In particular, in the theory of $p$-hypergroups, by a representation of the hypergroup $H$ on the Hilbert space $V$   is meant a symmetric homomorphism $f \mapsto \pi (f)$, which does not decrease the norm, of the Banach algebra $\mathscr{M}^b (H)$ of bounded measures to the algebra of bounded operators on $V$, such that the induced map $\mathscr{M}^b_+ S(V)$ is continuous (with respect to the weak operator topology on $S(V)$ and the topology on $\mathscr{M}^b_+ (H)$  described in point 5.3).   The matrix elements of this representation are bounded continuous functions on $H$ (cf., e.g., [136]).   The irreducible (symmetric) representations of an arbitrary p-hypergroup $H$ separate points in $H$ so that there are sufficiently many such representations [180].

The basic concepts, constructions, and results of the theory of unitary representations of groups generalize to the case of hypergroup algebras with involution; for example, analogs of the theorem on the connection of representations with positive definite functions, Bochner's theorem, Plancherel's theorem, etc. are valid.
\\

\underline{{\bf 6.2.}   Regular Representation.}     Let  $\mathscr{F}$ be a topological algebra, $\Phi$ be the dual space to $\mathscr{F}$, endowed with a locally convex topology compatible with the duality between $\Phi$ and $\mathscr{F}$, for example, the Mackey topology $\tau (\Phi , \mathscr{F} )$.    We denote by $R(f)$ the operator on $\Phi$ adjoint to the operator $g \mapsto gf$   and by $L(f)$ the operator adjoint to the operator $g \mapsto fg$.   The map $f \mapsto R(f)$ is a continuous representation of the algebra $\mathscr{F}$.   We call this representation {\it right-regular}.   We note that the correspondence $f \mapsto L(f)$  is an antirepresentation of the algebra $\mathscr{F}$.

If $\mathscr{F}$ coincides with the hypergroup algebra of the hypergroup $H$, while the g.t.o. $R^x$ and $L^x$ act on the space $\Phi$, then $R(f) = \int R^x df(x)$   and $L(f) = \int L^x df(x)$.     Hence the correspondence $x \mapsto R^x$ is a nondegenerate representation of the hypergroup $H$ and extends to a right-regular representation of the algebra $\mathscr{F}$.    The representation $x \mapsto R^x$  is called {\it the right-regular representation of the hypergroup $H$}.

If $H$ is a hypergroup with involution, then the correspondence $x \mapsto L^{x^*}$  is also a representation of the hypergroup $H$.   This representation is called {\it left-regular} and extends to the {\it left-regular representation} $f \mapsto L(f^*)$ of the algebra with involution $\mathscr{F}$.

Let us assume that the action of the operators $R^x, L^x, R(f), L(f)$ is defined not only on the space $\Phi$, but also in the Hilbert space $L^2 (H, m)$, where m is a right (left) invariant measure on $H$, cf. point 2.7.   In this case there arises on $L^2 (H, m)$ a symmetric representation $f \mapsto R(f)$  (respectively $f \mapsto L(f^*)$).    Representations of this type are also called regular.
\\

\underline{{\bf 6.3.}   Ideals and Representations.}   As in point 6.2 let $\mathscr{F}$ be a topological algebra, $\Phi$ be a locally convex space dual to $\mathscr{F}$.  By a {\it matrix element} of an arbitrary continuous representation $f \mapsto \pi (f)$ of the algebra $\mathscr{F}$ on the space $V$ is meant an element $\varphi \in \Phi$ such that  $\langle f, \varphi \rangle = \langle v', \pi (f)v \rangle$, where $v \in  V, v' \in  V'$.   The space $\Phi (\pi )$ in $\Phi$, spanned by matrix elements of the representation $\pi$, will be called the {\it space of matrix elements of this representation}.   We denote by  Ker $\pi$ the kernel of the representation $\pi$, i.e., the closed two-sided ideal $\{ f \in \mathscr{F} \colon \pi (f) = 0 \}$  in $\mathscr{F}$.   Continuous representations $\pi_1$ and $\pi_2$ of the algebra $\mathscr{F}$ are called {\it isomorphic} if Ker $\pi_1 =$ Ker$  \pi_2$.   For symmetric representations of $C^*$-algebras the concept of isomorphism was introduced by M. A. Naimark [69]; the general case was considered in [59-61].   For completely irreducible representations\symbolfootnote[2]{The representation $\pi \colon~ \mathscr{F} \to S(V)$ is completely irreducible, if its image is dense in $S(V)$; in this case $V$ does not contain nontrivial closed subspaces, invariant with respect to the operators $\pi (f)$, i.e., the representation $\pi$ is topologically irreducible.} isomorphism coincides with the {\it Fell equivalence} [123] (or {\it weak equivalence}), for irreducible finite-dimensional representations it consides with ordinary equivalence, defined by a similarity operator.   Conditions are given in [123] and [32] under which Fell equivalence (for irreducible representations) coincides with Naimark equivalence, which is popular in the theory of representations of semisimple groups.

We shall call a subspace \symbolfootnote[3]{By a subspace of a locally convex space we shall mean a closed linear subspace.} of $\Phi$, which is invariant with respect to all the operators $R(f)$ (respectively $L(f)$), {\it right-(left-)invariant}.   The next proposition follows from the definition of the operators $R(f)$ and $L(f)$.
\\

\underline{Proposition 5.}   The subspace $\Phi_1$ of $\Phi$ is right-(left-)invariant if and only if its orthogonal complement $\{f \in \mathscr{F} \colon \langle f, \varphi \rangle =0$ for all $\varphi \in \Phi_1 \}$ is a right (left) ideal in the algebra $\mathscr{F}$.
\\

Thus to the subrepresentations of the right-regular representation (i.e., to its restrictions to invariant subspaces of $\Phi$) correspond closed right ideals in $\mathscr{F}$ and by the Hahn-Banach theorem this correspondence is one-to-one.   It is clear that to maximal closed right ideals of $\mathscr{f}$ correspond topologically irreducible subrepresentations .
\\

\underline{COROLLARY 1.}   Closed two-sided ideals of $\mathscr{f}$ are in one-to-one correspondence with subspaces of $\Phi$, which are simultaneously right- and left-invariant.
\\

\underline{Proposition 6.}   The orthogonal complement in $\mathscr{F}$ of the space of matrix elements $\Phi (\pi)$ of the continuous representation $\pi$ coincides with the kernel of this representation Ker $\pi$.
\\

\underline{COROLLARY 2.}   Continuous representations of the algebra $\mathscr{f}$ are isomorphic if and only if their spaces of matrix elements coincide.
\\

Let $\delta$ be a fixed right identity in $\mathscr{F}$; we call the subalgebra $\widetilde{\mathscr{F}} =  \{f \in \mathscr{F} \colon~ \delta f = f \}$ the {\it basic subalgebra} of $\mathscr{F}$ and the subspace $\widetilde{\Phi} = \{ \varphi \in \Phi \colon~ L(\delta ) \varphi =\varphi \}$ the {\it basic subspace} of $\Phi$.  We call the representation $\pi \colon \mathscr{F} \to S(V)$  {\it non-degenerate}, if the operator $\pi (\delta )$ is the identity operator.   This terminology agrees with the definitions introduced above for hypergroup algebras.   The algebra $\tilde{\mathscr{F}}$  has a two-sided identity $\delta$ and is isomorphic to the quotient-algebra $\mathscr{F} / J$,  where $J$ is the orthogonal complement of $\widetilde{\Phi}$.   The nondegenerate representations of the algebra $\mathscr{F}$  are in one-to-one correspondence with representations of the algebra $\tilde{\mathscr{F}}$,   as in Theorem 3.   The closed ideals (right, left, two-sided) of  $\tilde{\mathscr{F}}$  are in one-to-one correspondence with the invariant subspaces of $\widetilde{\Phi}$ (respectively, right-, left-, or bilaterally invariant).

In particular, if $\mathscr{F} = \mathscr{M} (H)$  where $H$ is a locally compact hypergroup, then the space of matrix elements of any nondegenerate representation of $H$ is contained in the basic subspace $\widetilde{C} (H)$, which coincides with the space of matrix elements of the right-regular representation.   We note that for an arbitrary hypergroup it is impossible to define the tensor product of representations.   Hence the collection of all matrix elements will not necessarily be an algebra with respect to multiplication of functions.

In what follows, when we speak of representations of hypergroups having certain properties or others (for example, being irreducible, equivalent, etc.), we shall have in mind that the corresponding representations of hypergroup algebras have these properties.
\\

\underline{{\bf 6.4.}   Characters.}   Reduction to hypergroup algebras lets one carry over the concept of the character of a representation to the case of hypergroups, where it is treated as a linear functional on an ideal in the hypergroup algebra $\mathscr{F}$ (or on an  $\mathscr{F}$-bimodule).   If the character $\chi$ of a topologically irreducible representation $\pi$ is defined on a dense ideal in $\mathscr{F}$ and $\chi (f)$ coincides with the trace of a nuclear operator $\pi (f)$, then the representation $\pi$ is completely irreducible and is determined by its character $\chi$ up to Fell equivalence [63].   Cf. [60] for a more general approach and other conditions under which a representation is determined by its character up to isomorphism.  As in the case of Lie groups, under certain conditions the character can be considered as a generalized function on a hypergroup.

In what follows, by a {\it character of a hypergroup} $H$ we shall mean a nonzero one-dimensional representation (which can be identified with its character in the sense considered above or with a matrix element).   Thus, a character can be considered as a function $\chi (x)$ on $H$, satisfying the equation
$$
R^x \chi (x) = \chi (y) \chi (x). \eqno(6.1)
$$
If $H$ is a hypergroup with involution and the representation  $x \mapsto \chi (x)$  is symmetric, then $\chi (x^*) = \bar{\chi} (x)$; in this case the character $\chi$ is called {\it Hermitian}.   Characters of commutative hypergroups were studied by many authors, cf., e.g., [3, 7, 41, 48-51, 55, 56, 61, 112, 117, 136, 169].
\\

\underline{Example.   Spherical Functions.}   Let $H$ be the hypergroup generated by the Gelfand pair $(G, K)$, cf. above, points 2.6 and 5.4.   The characters of this hypergroup correspond to continuous functions $\varphi (x)$ on $G$, satisfying the relation
$$
\int_K \varphi (xky) dk = \varphi (x) \varphi (y) \eqno(6.2)
$$
for all $x, y \in K$ and the normalized Haar measure dk.   (6.2), which follows from (2.4) and (6.1), means that $\varphi (x)$ is a spherical function on $G$, cf., e.g., [23, 10, 69, 84, 85].
\\

\underline{{\bf 6.5.}   Generalized Fourier Transform.}   Let $H$ be a hypergroup with a hypergroup algebra $\mathscr{F}$ and suppose that with each element $\zeta$ of some set $Z$ there is associated a representation $\pi_{\zeta}$ of the hypergroup $H$.   By the {\it generalized Fourier transform} we mean the linear operator, which carries the functional $f \in \mathscr{F}$  into the function  $f(\zeta ) = \pi_{\zeta} (f)$   assuming operator values.   It is clear that here the generalized convolution goes into the product of the operator-valued functions.   Such a construction is standard in the representation theory.

If the set $Z$ consists of characters of the hypergroup $H$, then $f(\zeta )$ is a scalar function and the generalized convolution goes into ordinary multiplication of functions.  Here the basic subalgebra $\widetilde{\mathscr{F}}$ is mapped isomorphic-ally onto some algebra of functions if and only if linear combinations of characters from $Z$ are dense in the basic subspace $\widetilde{\Phi}$, where $\Phi = \mathscr{F}'$.  For example, if with each complex number $\zeta$ one associates the character $x \mapsto e^{ix \zeta}$ of the group ${\bf R}$ of real numbers, then for the algebra $\mathscr{D} ({\bf R})$  the construction indicated leads to the ordinary Fourier-Laplace transform of generalized functions with compact support.   If the hypergroup $H$ is generated by a Gelfand pair and $Z$ consists of Hermitian characters of this hypergroup, then the generalized Fourier transform reduces to the spherical Fourier transform, cf., e.g., [7, 10, 84, 85, 93, 131, 136, 146, 169].   Other examples are given below.
\\

\underline{{\bf 6.6.}   Spectral Analysis and Synthesis of Ideals and Representations.}  It is clear that the description of continuous representations of the topological algebra $\mathscr{F}$  up to isomorphism reduces to the description of closed two-sided ideals of $\mathscr{F}$. In what follows we shall consider only continuous representations and closed two-sided ideals.   We are particularly interested in (nontrivial) maximal, primitive and primary ideals.   {\it Primitive} ideals are kernels of irreducible representations, and an ideal which is contained in a unique maximal ideal is called {\it primary} (we note that this definition, which is customary in the theory of topological algebras, does not agree with the definition of a primary ideal in abstract algebra; sometimes an ideal which is contained in a unique primitive ideal is called a primary ideal in a topological algebra).   A representation $\pi$ is called {\it prime} (respectively {\it primary}), if its kernel is a maximal (primary) ideal.   A finite-dimensional representation is prime if and only if it is isomorphic to an irreducible finite-dimensional representation; cf. [60] on prime representations.

We shall say that the representation $\pi_1$ is {\it contained in the representation} $\pi_2$, if Ker $\pi_1 \supset$ Ker $\pi_2$.   In this case the space of matrix elements of the representation $\pi_2$ contains the space of matrix elements of the representation $\pi_1$.   If $\pi_1$ is a subrepresentation of $\pi_2$, then clearly $\pi_1$ is contained in $\pi_2$; the converse assertion is true only "up to isomorphism."   It is clear that a representation is primary if it contains a unique prime representation.

We fix some collection $\Gamma = \{ I_{\gamma} \}$ of ideals in $\mathscr{F}$ and for any ideal $I_{\gamma} \in \Gamma$ we denote by $\pi_{\gamma}$ a representation such that $I_{\gamma} =$ Ker $\pi_{\gamma}$.   It is appropriate to choose $\Gamma$ so that this collection is visible, consists of ideals of relatively simple structure, and other ideals can be described conveniently as intersections of ideals of $\Gamma$ in the spirit of the famous theorem of Lasker-Noether.   We shall say that the ideal $I$ admits {\it spectral synthesis} (with respect to $\Gamma$), if $I$ coincides with the intersection of those ideals of $\Gamma$ which contain it; a {\it representation $\pi$ admits spectral synthesis} if the ideal Ker $\pi$ admits it.   The next proposition follows from Propositions 5 and 6 and the Hahn-Banach theorem.
\\

\underline{Proposition 7.}   The representation $\pi$ admits spectral synthesis if and only if any matrix element of this representation can be approximated in $\mathscr{F}'$ by matrix elements   of those representations $\pi_{\gamma}$  which are contained in $\pi$.
\\

For example, for representations of a locally compact hypergroup one is concerned with approximation of matrix elements in the topology of uniform convergence on compact subsets of $H$.

It is natural to treat the problem of {\it spectral analysis} as the problem of describing ideals of $\Gamma$, containing a given ideal, or representations $\pi_{\gamma}$ contained in a given representation.   The scheme described of one of the versions of harmonic analysis, connected with the Delsarte-Schwartz theory of mean-periodic functions and the theory of differential equations with constant coefficients and convolution equations, allows one to consider results which at first glance have little in common between them, from a single point of view.   A vast literature is devoted to other questions of harmonic analysis on hypergroups (the theory of almost periodic functions, expansions in orthogonal systems of functions, etc., cf., points 1.6, 5.3, 6.1, 6.8).
\\

\underline{{\bf 6.7.}   Examples of Spectral Synthesis.}   The theory of representations of compact groups gives the classical example of spectral synthesis.   If $G$ is a compact group, then in the topological group algebra $\mathscr{M} (G)$ consisting of measures with compact support, any primary ideal is a maximal ideal of finite codimension.   The Peter-Weyl theorem is equivalent to the assertion that any ideal\symbolfootnote[1]{We recall that only closed two-sided ideals are considered.} in $\mathscr{M} (G)$ is the intersection of the maximal ideals containing it.

A generalization of the Peter-Weyl theorem to the case of compact hypergroups is given by B. M. Levitan (cf. [55, 56]; the Peter-Weyl theory and harmonic analysis on p-hypergroups were considered in [136, 180]; cf. also [3, 8, 7, 14]).   If there is given on the compact hypergroup $H$ a finite measure, satisfying certain additional conditions, then any continuous function from the basic subspace $\widetilde{C} (H)$ can be approximated in the topology of uniform convergence by matrix elements of finite-dimensional irreducible representations.   One can deduce from this that any ideal in the basic subalgebra $\tilde{\mathscr{M}} (H)$   is the intersection of the maximal ideals containing it. Analogous results are also valid for other hypergroup algebras, for example, for the Hilbert algebra $L^2 (H)$.

For the group $G$ of conformal transformations of the unit disc, which is locally isomorphic to the group $SL(2, R)$, Ehrenpreis and Mautner [119] described the ideals in the group algebra  $\mathscr{D} (G)$   (we note that for any Lie group $G$ the ideals in  $\mathscr{D} (G)$   are in one-to-one correspondence with the ideals in $\mathscr{M} (G)$).   It turns out that there are ideals which are not intersections of primary ideals, and sufficient conditions are given in [119] under which an ideal admits spectral synthesis.   Using the concept of g.t.o., P. K. Rashevskii found exhaustive results on spectral analysis and synthesis for other function spaces\symbolfootnote[2]{Here instead of group and hypergroup algebras, bimodules over them were considered.} on $SL(2, R)$.   Cf. [76, 87] for the further development of problems of this kind for other semisimple groups of rank 1.   For groups of arbitrary rank additional difficulties arise, which are familiar from the "multidimensional" theory of mean-periodic functions.

We consider the topological algebras $\mathscr{M} ({\bf R}^n ), \mathscr{D} ({\bf R}^n ), \mathscr{A} ({\bf C}^n)$ connected with ordinary translations in an n-dimensional vector space.  The multiplication in these algebras coincides with the usual convolution of the corresponding measures, generalized functions, or analytic functionals.   The problem of describing the ideals in these topological algebras and of the spectral synthesis for such ideals is equivalent with the basic problem of the theory of mean-periodic functions.   For $n = 1$, this problem is solved in [165].   For any complex number $\zeta$ and integer $k \geq 0$ we denote by $C(\zeta , k)$ the subspace of $C({\bf R})$   spanned by the quasimonomials $e^{\zeta t} , te^{\zeta t} , \ldots , t^k e^{\zeta t}$ and by $I(\zeta , k)$ we denote the ideal in $\mathscr{M} ({\bf R})$  orthogonal to $C(\zeta , k)$.   Ideals of this type (and only these) are primary.   Any nontrivial ideal in $\mathscr{M} ({\bf R})$  is the intersection of a no more than countable family of primary ideals, and each nontrivial closed subspace of $C({\bf R})$, which is invariant with respect to translations, is spanned   by a no more than countable set of quasimonomials.   Analogous results are valid for the algebras $\mathscr{D} ({\bf R})$   and $\mathscr{A} ({\bf C})$.     For $n > 1$ the situation is more complex.   Schwartz assumed that in the multidimensional case also the invariant subspaces are spanned by quasimonomials.   In this case the corresponding ideals admit spectral synthesis.   Although a counterexample has been constructed to Schwartz' conjecture [26], this conjecture is proved for many important special cases, cf., e.g., [148, 116, 118, 75, 70].   For spaces of solutions of homogeneous systems of linear partial differential equations with constant coefficients, stronger results have been found which can be interpreted as the decomposition of nonunitary representations into the direct integral of primary components [75, 118].   Cf. [44-46] for generalizations of these results, for example, to the case of symmetric spaces.   For certain one-dimensional hypergroups, the basic results of the theory of mean-periodic functions (including the primary synthesis) are found in [41].   Analogous results are found in [91] for the hypergroup $H = K\smallsetminus G/K$ where $G$ is a semisimple linear Lie group and $K$ is a maximal compact subgroup of it.

The following example is described in [61].   Let the real line ${\bf R}$ be endowed with the structure of the hypergroup described in point 2.4, and let $\mathscr{M}$ be the corresponding topological hypergroup algebra, consisting of measures with compact support.   The basic subspace $\widetilde{C} ({\bf R})$ consists of all even continuous functions, and the basic subalgebra $\tilde{\mathscr{M}}$ consists of all even measures with compact support.   Equation (6.1), which defines the characters, has, in the present case, the form $\chi (x+y) + \chi (x-y) = 2\chi (x) \chi (y)$ and the characters are the functions $\chi_{\zeta} (x) = cos (\zeta x)$  where $\zeta$ is any complex number.   The generalized Fourier transform $f \mapsto f(\zeta ) = \int \cos (\zeta x ) df(x)$  maps $\tilde{\mathscr{M}}$ one-to-one to the algebra consisting of even entire analytic functions of exponential growth (not all).   Any primary ideal in $\tilde{\mathscr{M}}$ has the form $\{ f \in \tilde{\mathscr{M}} \colon f (\zeta ) = 0, f' (\zeta ) = 0, \ldots , f^{(k)} (\zeta ) = 0 \}$ where the complex number $\zeta$ and the integer $k \geq 0$ are fixed.   Any nontrivial ideal in $\tilde{\mathscr{M}}$  is the intersection of a no more than countable collection of primary ideals.   Any nontrivial subspace of $\widetilde{C} ({\bf R})$ which is invariant with respect to g.t.o., is spanned by a no more than a countable set of functions of the form $\cos (\zeta x ), x^{2n - 1} \sin (\zeta x), x^{2n} \cos (\zeta x )$ where $\zeta \in {\bf C}, n= 1, 2, \ldots$.   Analogous results are valid for smooth and holomorphic functions.

The more general case of g.t.o. with Sturm -- Liouville generator is considered in [27], where examples of nonprimary synthesis are given.   One can extract additional examples, for example from [59, 61, 83, 109].
\\

\underline{{\bf 6.8.}   Plancherel's Theorem and the Inversion Formula.}   Suppose given on certain spaces $H$ and $\widehat{H}$ positive measures $m$ and $\mu$ respectively.   Let us assume that with the help of the function $u(x, \chi )$ defined on $H \times \widehat{H}$, there is given an "abstract Fourier transform"
$$
\varphi (x) \mapsto \hat{\varphi} (\chi ) = \int \varphi (x) \overline{u(x, \chi )} dm(x), \eqno(6.3)
$$
which defines an isomorphism of the Hilbert spaces $L^2 (H, m)$ and $L^2 (\widehat{H} , \mu )$, so that the "Plancherei formula" is valid:
$$
\int \varphi (x) \overline{\psi (x)} dm(x) = \int \hat{\varphi} (\chi) \overline{\hat{\psi} (\chi )} d\mu (\chi ) .\eqno(6.4)
$$
Let us assume in addition that the {\it inversion formula}
$$
\varphi (x) = \int \hat{\varphi} (\chi ) u (x, \chi ) d\mu (\chi) \eqno(6.5)
$$
is valid, which one can get by formal substituting the delta-function for $\psi$ in (6.4) and considering the form of the transformation (6.3).   If the measure $\mu$ is discrete, then (6.5) gives the expansion of the function $\varphi$ in an "abstract Fourier series."

It turns out (cf. [51, 55, 56]) that $H$ has the structure of a commutative hypergroup, if for some point $e \in H$ and for all $\chi \in \widehat{H}$ one has $u(e, \chi ) = 1$ (the element $e$ is the identity in $H$; sometimes one is naturally restricted by the requirement of the existence of an approximate identity).   Here the g.t.o. are defined by
$$
R^y \varphi (x) = \int \hat{\varphi} (\chi ) u(y, \chi ) u(x, \chi ) d\mu (\chi ) .
$$
Analogously one defines a hypergroup structure on $H$.   Hence hypergroups and g.t.o. arise naturally in problems connected with expansion in orthogonal systems of functions, the spectral theory of operators, etc., which guarantees a large circle of applications of the theory of hypergroups and g.t.o.

For a large class of hypergroups with invariant measure the transformation (6.3) can be defined as a {\it generalized Fourier} transform (cf. point 6.5), applied to the measure $\varphi (x)m$, to give explicitly the space $\widehat{H}$  and the measure $\mu$ (called the {\it Plancherel measure})   and to prove the Plancherel formula and the inversion formula. In the commutative case $\widehat{H}$ consists of Hermitian characters (or of one-dimensional representations of the corresponding hypergroup algebra with involution) and $u(x, \chi ) = \chi (x)$, where $\chi$ is a character from $\widehat{H}$.   The Plancherel theorem and inversion formula for commutative g.t.o. were first proved by B. M. Levitan [47-50], cf. also [3, 7, 136, 169].   In [17] a generalization of the Plancherel theorem and inversion formula to a certain class of (generally) noncommutative hypergroups is given\symbolfootnote[1]{The class of commutative hypergroups for which one is able to prove an "abstract" Plancherel theorem is extended simultaneously.}.   This class includes, for example, all p-hypergroups with invariant measure and all locally compact groups.   Thus, this result generalizes the theorem proved by Segal [166] and Tatsuuma [173] for locally compact groups.

We note that the existence of the inversion formula is a distinctive property of hypergroups, while the Plancherel theorem is also valid for objects of a more general nature, connected with Hilbert algebras, for example, for groupoids with measure, cf. [159].

\section*{\small 7.     Conclusion}
The theory of hypergroups is a relatively new domain of mathematics, comparable, in breadth of content of its material, interest, and richness of applications, with such classical domains as the theory of Lie groups or the theory of Hilbert spaces, but considerably less developed.   The Lie theory for hypergroups (especially in relation to procedures for constructing global objects from given infinitesimal objects, for example, commutation relations) is far from to be complete.   Comparatively little concrete material has been accumulated in the theory of representations of noncommutative hypergroups.   In particular, there is interest in the calculation of concrete Plancherel measures and their supports, and in "nonunitary" theory, in description   of collections of ideals in hypergroup algebras and invariant subspaces in function spaces on hypergroups, guaranteeing the meaningfulness of spectral synthesis.

One can hope that the contemporary theory of hypergroups and g.t.o. will enter as a fragment or the skeleton of a broader and more general theory, which will include, for example, the theory of supergroups also. Generalizations of the concept of hypergroup arise naturally in Lie theory, duality theory, group representation theory.   The direction of further development of the theory of hypergroups will apparently be determined to a considerable degree by impulses connected with the applications of this theory to mathematical physics.

{\it Grigory L. Litvinov, Moscow

E-mail: glitvinov@gmail.com}

\end{document}